\documentclass[12pt,letter]{article}  
\usepackage[T2A]{fontenc}
\usepackage[cp866]{inputenc}
\usepackage[russian,english]{babel}
\usepackage[tbtags]{amsmath}
\usepackage{amsfonts,amssymb,mathrsfs,amscd,comment}
\usepackage{amsthm,graphicx} 
\usepackage{hyperref}
\def\NoBlackBoxes{\overfullrule0pt}

\NoBlackBoxes
\theoremstyle{plain}
\newtheorem{theorem}{Theorem}
\newtheorem{theoremML}{Theorem (Mittag-Leffler)}

\newtheorem{lemma}{Lemma}

\theoremstyle{definition}
\newtheorem{definition}{Definition}

\newtheorem{remark}{Remark}
\newtheorem{example}{Example}

\theoremstyle{main}
\def\bad{\spaceskip=0.33emplus0.6emminus0.15em\immediate\write5{\string\bad}}
\let\leq\leqslant
\def\ZZ{\mathbb Z}
\def\({\left(}
\def\){\right)}
\def\[{\left[}
\def\]{\right]}
\def\<{\left\langle}
\def\>{\right\rangle}
\def\Re{\operatorname{Re}}

\def\mdet{\operatorname{det}}
\def\Sh{\operatorname{Sh}}

\def\mdeg{\operatorname{deg}}
\def\mcap{\operatorname{cap}}
\def\supp{\operatorname{supp}}
\def\const{\operatorname{const}}

\def\NN{\mathbb N}
\def\RR{\mathbb R}
\def\CC{\mathbb C}
\def\PP{\mathbb P}
\def\DD{\mathbb D}
\def\QQ{\mathbb Q}
\def\KK{\mathbf K}
\def\bU{\mathbf U}
\def\SZ{\mathscr Z}

\def\HH{\mathscr H}
\def\FF{\mathbf F}
\def\bM{\mathbf M}
\def\EE{\mathbf E}

\def\MM{\mathscr M}

\def\RS{\mathfrak R}

\def\zz{\mathbf z}

\def\pp{\mathbf p}
\def\aa{\mathbf a}
\def\bb{\mathbf b}

\let\pfi\varphi
\let\leq\leqslant
\let\ge\geqslant
\let\myh\widehat\let\myt\widetilde\let\myo\overline
\def\bad{\spaceskip=0.33emplus0.6emminus0.15em\immediate\write5{\string\bad}}

\def\nn{\mathbf n}
\def\bGamma{\boldsymbol\Gamma}
\def\myk{\mathfrak K}
\begin{document}

\selectlanguage{english}

\title{Hermite--Pad\'e polynomials and analytic continuation:
new approach and some results}

\author{Sergey~P.~Suetin}


\maketitle

\markright{Analytic continuation}


\begin{abstract}
We discuss a new approach to realization of the well-known Weierstrass's programme on
efficient continuation of an analytic element
corresponding to a~multivalued analytic function with finite number of branch points.
Our approach is based on the use of Hermite--Pad\'e polynomials.

Bibliography:~\cite{Sze21}~titles.

\end{abstract}

\footnotetext[0]
{Steklov Mathematical Institute of the Russian Academy of Sciences, Moscow,
suetin@mi.ras.ru}

\footnotetext[0]{
This research was carried out with partial financial support of the Russian Foundation for Basic
Research (grant no.~18-01-00764).}

\section{Introduction and statement of the problem}\label{s1}

\subsection{}\label{s1s1}
The problem of analytic continuation of power series beyond its disc of convergence
is a~well-known problem in complex analysis, which also has
important applications. In the present paper, we discuss
this problem from the point of view of Weierstrass's approach to the concept
of analytic functions, which depends on the local representation of
an analytic function as a~power series with centre at some point of the Riemann sphere~$\myo\CC$.
For more on this approach, see, \textit{prima facie}, the book \cite{Bib67},
and also \cite{Bor05},~\cite{Pai05}, \cite[Ch.~8]{Mar68},
\cite{Ara84},~\cite{ArLu02},~\cite{Ara03}.

The purpose of the present paper is to introduce and briefly discuss a~new
approach to the efficient solution of the analytic continuation problem for power series
residing on the use of Hermite--Pad\'e polynomials and some extensions thereof.
In \S\S~\ref{s2} and~\ref{s3} we formulate some new theoretical results in this direction
which we have been able to derive so far. In \S~\ref{s5}  we briefly discuss
the most straightforward, to our opinion, applications of our theoretical results in the study of
applied problems.
In particular, in~\S~\ref{s5} we illustrate our approach and theoretical results
on some numerical examples pertaining to the model class of multivalued analytic functions
considered in the present paper (see \S\,\ref{s1s2} below).

We plan to prove the theoretical results of the present study
in the separate papers \cite{IkKoSu18} and~\cite{LoKoSu18}.
We also intend to
give a~more detailed analysis of possible applications of the method proposed here in applied problems
involving analytic continuation.

The new approach proposed below to efficient solution of the analytic continuation problem for power series
beyond its disc of convergence
will be demonstrated (from the theoretical as well as the numerical standpoints)
on an example of some ``model'' class of multivalued analytic functions based on the use
of the inverse of the Zhukovskii function (see representation \eqref{1} below).
This class was first introduced in \cite{Sue17} (see also the papers \cite{Sue17b},~\cite{Sue18}), where
it was denoted by~$\SZ$. Below, this notation will be retained, but when required
the parameters in representation \eqref{1} will be refined.

\subsection{}\label{s1s2}
Let $\Delta_1:=[-1,1]$, $\pfi(z):=z+(z^2-1)^{1/2}$, $z\notin \Delta_1$, be
the inverse of the Zhukovskii function; here and in what follows we
choose the branch of the root function such that
$(z^2-1)^{1/2}/z\to1$ as $z\to\infty$. The function $\pfi(z)$ is meromorphic
in the domain $D_1:=\myo\CC\setminus\Delta_1$ ($\pfi$ maps $D_1$ conformally and univalently onto
the exterior  $\myo\CC\setminus\myo\DD$ of the unit disc
$\DD:=\{z:|z|<1\}$). In the domain $\myo\CC\setminus\{-1,1\}$ the function $\pfi$
is already a~multivalued (more precisely, two-valued) analytic function.

Let $A_j\in\CC$, $j=1,\dots,m$, be arbitrary pairwise distinct complex numbers such that $|A_j|>1$.
Consider the function
\begin{equation}
f(z):=\prod_{j=1}^m\(A_j-\frac1{\pfi(z)}\)^{\alpha_j},
\quad z\in D_1,
\label{1}
\end{equation}
where $\alpha_j\in\CC\setminus\ZZ$, $\sum_{j=1}^m\alpha_j\in\ZZ$ (note that for the above
branch of the root function  we have $1/\pfi(z)=z-(z^2-1)^{1/2}$). Since
$(A_j-1/\pfi(z))\neq0$ for $z\in D_1$ under the above conditions on $A_j$ and
with the above choice of the branch of the root function, we can
find a~holomorphic  (i.e., singlevalued analytic) branch
$f\in\HH(D_1)$ of the function~$f$ in the domain $D_1$.
In the domain
$\myo\CC\setminus\{-1,1\}$ the analytic function~$f$ is not anymore singlevalued.
In addition to the second-order branch points $z=\pm1$, this function has branches at the points
$a_j=(A_j+1/A_j)/2\in\CC\setminus\Delta_1$, $j=1,\dots,m$, of (in general) infinite order (if
$\alpha_j\in\CC\setminus\QQ$ for the corresponding~$A_j$). Thus, the total
set of branch points of a~function~$f$ of
the form~\eqref{1} is $\Sigma=\Sigma(f)=\{\pm1,a_j,j=1,\dots,m\}$.

In what follows, for an arbitrary $n\in\NN$, we denote by $\PP_n:=\CC[z]$ the space of all
algebraic polynomials of degree $\leq{n}$,
$\PP^{*}_n:=\PP_n\setminus\{0\}$. Given an arbitrary polynomial $Q\in\CC[z]$,
$Q\not\equiv0$, we let
$$
\chi(Q):=\sum_{\zeta:Q(\zeta)=0}\delta_\zeta
$$
demote the measure counting the zeros (with multiplicities) of the polynomial~$Q$.

The weak convergence in the space of measures will be denoted by
``$\overset{*}\longrightarrow$''.

For an arbitrary (positive Borel) measure $\mu$,
$\supp{\mu}\subset\CC$, we denote by~$V^\mu$ the logarithmic potential of~$\mu$,
$$
V^\mu(z):=\int\log\frac1{|z-\zeta|}\,d\mu(\zeta).
$$

By  ``$\overset\mcap\longrightarrow$'' we shall denote the convergence with respect to the
(logarithmic) capacity on compact subsets of some domain.

Below, ``$\sqrt{\,\cdot\,}\,$'' means the principle square root
of a~nonnegative real number, $\sqrt{a^2}=|a|$, $a\in\RR$.

Given a~function $f\in\HH(D_1)$, we fix some germ\footnote{Throughout, the
terms
``germ'', ``analytic element'' and ``power series'' are used to denote the same object: a~convergent
power series with centre at some point of the extended complex plane $\myo\CC$ and nonzero convergence radius.}
$$
f_\infty=\sum_{k=0}^\infty\frac{c_k}{z^k},
$$
holomorphic at the point $z=\infty$, $f_\infty\in\HH(\infty)$.

For a fixed germ $f_\infty\in\HH(\infty)$ and an arbitrary $n\in\NN$, we denote by
$P_n$ and $Q_n$, $P_n,Q_n\in\PP_n:=\CC[z]$, $Q_n\not\equiv0$, the polynomials defined
(not uniquely) from the relation
\begin{equation}
(Q_nf_\infty-P_n)(z)=O\(\frac1{z^{n+1}}\),\quad z\to\infty.
\label{2}
\end{equation}
The rational function $[n/n]_{f_\infty}:=P_n/Q_n$, which is uniquely determined,
is called the {\it diagonal Pad\'e approximant} (the ``diagonal PA'' or simply ``PA'')
of the germ $f_\infty$ at the point $z=\infty$.

The next facts \eqref{3}--\eqref{5}  follow from
Stahl's theory (see \cite{Sta97b},~\cite{Sta12}, and also \cite{ApBuMaSu11}
and~\cite{Sue15}):
\begin{equation}
\frac1n\chi(P_n),\frac1n\chi(Q_n)\overset{*}\longrightarrow\tau_{\Delta_1},
\quad n\to\infty,
\label{3}
\end{equation}
here $\tau^{\vphantom{p}}_{\Delta_1}=\dfrac1\pi\dfrac{dx}{\sqrt{1-x^2}}$ is the unit Robin
measure of the interval $\Delta_1$; i.e., $V^{\tau^{\vphantom{p}}_{\Delta_1}}(x)\equiv\const=\gamma_{\Delta_1}$,
and $\gamma_{\Delta_1}$ is the Robin constant for the interval $\Delta_1$;
\begin{equation}
[n/n]_{f_\infty}(z)\overset\mcap\longrightarrow f(z)\in\HH(D_1),
\quad n\to\infty,\quad z\in D_1;
\label{4}
\end{equation}
the convergence in~\eqref{4} can be characterized by the relation
\begin{equation}
\bigl|f(z)-[n/n]_{f_\infty}(z)\bigr|^{1/n}
\overset\mcap\longrightarrow e^{-2g_{\Delta_1}(z,\infty)}<1,
\quad n\to\infty,\quad z\in D_1,
\label{5}
\end{equation}
where $g_{\Delta_1}(z,\infty)$ is the Green function of the domain $D_1$ with singularity
at the point at infinity  $z=\infty$. Relation \eqref{5} means that the PA
$[n/n]_{f_\infty}$ converges to the function $f\in\HH(D_1)$
with the rate of a~geometric progression with ratio $\delta(z):=e^{-2g_{\Delta_1}(z,\infty)}<1$.

So, from \eqref{3}--\eqref{5} it follows that the diagonal PA, which are
constructed only from {\it local} data (the germ $f_\infty$ of the function~$f$ at the point $z=\infty$),
recover the function~$f$ in the entire domain $D_1$. From the viewpoint of
Stahl's theory, the domain $D_1=D_{\max}(f_\infty)$ is a~``maximal'' domain
in which the diagonal PA converge (in capacity).
The ``maximality'' of~$D_1$ is understood in the sense that the boundary $\partial
D_1=\Delta_1$ has minimal capacity among the boundaries of all admissible
domains~$G$ for the germ $f_\infty$; i.e.,
$$
\mcap(\partial D_1)=\min\{\mcap(\partial G): G\ni\infty, f_\infty\in\HH(G)\}.
$$
From~\eqref{3} it follows that ``almost all'' zeros and
poles\footnote{More precisely, all but at most  $o(n)$ many zeros and poles.}
of the diagonal PA $[n/n]_{f_\infty}$ accumulate to the closed interval
$\Delta_1=\partial D_1$.
In the above context, the interval $\Delta_1$ is called
the {\it  Stahl compact set}, and the domain $D_1$, the {\it Stahl domain}. Of course,
in the framework of the general Stahl theory, which holds for an arbitrary multivalued analytic function
with finite number of branch points, these two concepts are more meaningful than in the very particular
case considered here.

Thus, the example of this particular class $\SZ$ shows that the diagonal PA, as
constructed only from local data, are capable of solving the following problems:

1) recover the Stahl compact set, and therefore, the Stahl domain
$D_1=\myo\CC\setminus\Delta_1$ at which $f_\infty\in\HH(D_1)$,
from limit distribution of its zeros and poles;

2) provide an efficient (singlevalued) analytic continuation of a~given germ
$f_\infty\in\HH(\infty)$ to the domain $D_1$ as a~holomorphic function $f(z)$, $z\in
D_1$.

\smallskip

It follows from what has been said that for all functions~$f$ from the class~$\SZ$
the Stahl compact set $S(f_\infty)$ is equal to $\Delta_1$ (i.e., it is
independent\footnote{For the already chosen
branch of the inverse of the Zhukovskii function $\pfi(z)$.}
of the function $f\in\SZ$) and that the so-called
``active'' branch points of the germ $f_\infty$ are the points $z=\pm1$ (for the definition of
active branch points, see \cite{Sta12}).
Note that for $\alpha_j\in\CC\setminus\QQ$ the Riemann surface (or, simply, R.s.) $\RS=\RS(f)$ of a~function
$f\in\SZ$ is infinitely-sheeted. The Stahl domain
$D_1$ can be looked upon as the ``first'' sheet of this R.s.\ $\RS(f)$. Hence
the result of the application of Stahl's theory to the germ $f_\infty$ of the function $f\in\SZ$
can be interpreted as follows. The diagonal PA $[n/n]_{f_\infty}$, as constructed from
the germ~$f_\infty$, recover the function~$f$ only on the first sheet of the
R.s.\ $\RS(f)$, while all inactive branch points $\{a_j,j=1,\dots,m\}$ of the
function~$f$ lie on the ``other'' sheets of this R.s.\ $\RS(f)$.
Therefore, in contrast to the points $z=\pm1$, the remaining branch points $\{a_j,j=1,\dots,m\}$
turn out to be inaccessible in the sense of their recovery from zeros and poles
of the diagonal Pad\'e approximants.

\smallskip

This suggests the following fairly natural questions.

1) {\it Is it possible, from a~given germ $f_\infty\in\HH(\infty)$ of a~multivalued
analytic function~$f$ (in particular, from the class~$\SZ$), to recover the remaining branch points of~$f$
which are inactive for the germ $f_\infty$ in the framework of Stahl's theory}?

2) {\it Is it possible, from a~given germ $f_\infty$, to recover the function~$f$ on
``other sheets'' of its R.s.\ $\RS(f)$, rather than only on the first sheet,
as is inherent in Stahl's theory}?

\smallskip

Note that the above questions fit well the lines of the general {\it
Weierstrass programme} (see \cite{Ara84},~\cite{ArLu02},~\cite{Ara03} and the references
given there), which is aimed at extracting all properties of the so-called
{\it global} analytic function directly in terms of its specific germ (i.e., in the actual fact, in
terms of the corresponding Taylor coefficients of power series).

\smallskip

In the present paper we consider a~very special class of multivalued
analytic functions: the analytic functions given by the explicit representation~\eqref{1}.
Nevertheless, this class is fairly representative in the sense that it shows exactly which
advantages come from the  use of rational approximants based on Hermite--Pad\'e polynomials in comparison with
Pad\'e diagonal approximants. Namely, below we shall formulate some theoretical results from which
it follows that at least in the class~$\SZ$ the above questions 1)--2) can be answered affirmatively.

\section{The real case}\label{s2}

\subsection{}\label{s2s1}
Assume that in \eqref{1}\enskip  $m=2$, $\alpha_1=\alpha_2=-1/2$, and $A_1,A_2$ are real numbers
such that $1<A_1<A_2$. In what follows, it is convenient to put $A_1=A$,
$A_2=B$. Hence $1<A<B$, and
\begin{equation}
f(z)=\[\(A-\frac1{\pfi(z)}\)\(B-\frac1{\pfi(z)}\)\]^{-1/2},
\quad z\in D_1,
\label{6}
\end{equation}
$\Sigma(f)=\{\pm1,a,b\}$, where $a=(A+1/A)/2$, $b=(B+1/B)/2$, $1<a<b$.
We set $\Delta_2=[a,b]$.

\begin{figure}[!ht]
\centerline{
\includegraphics[width=15cm,height=15cm]{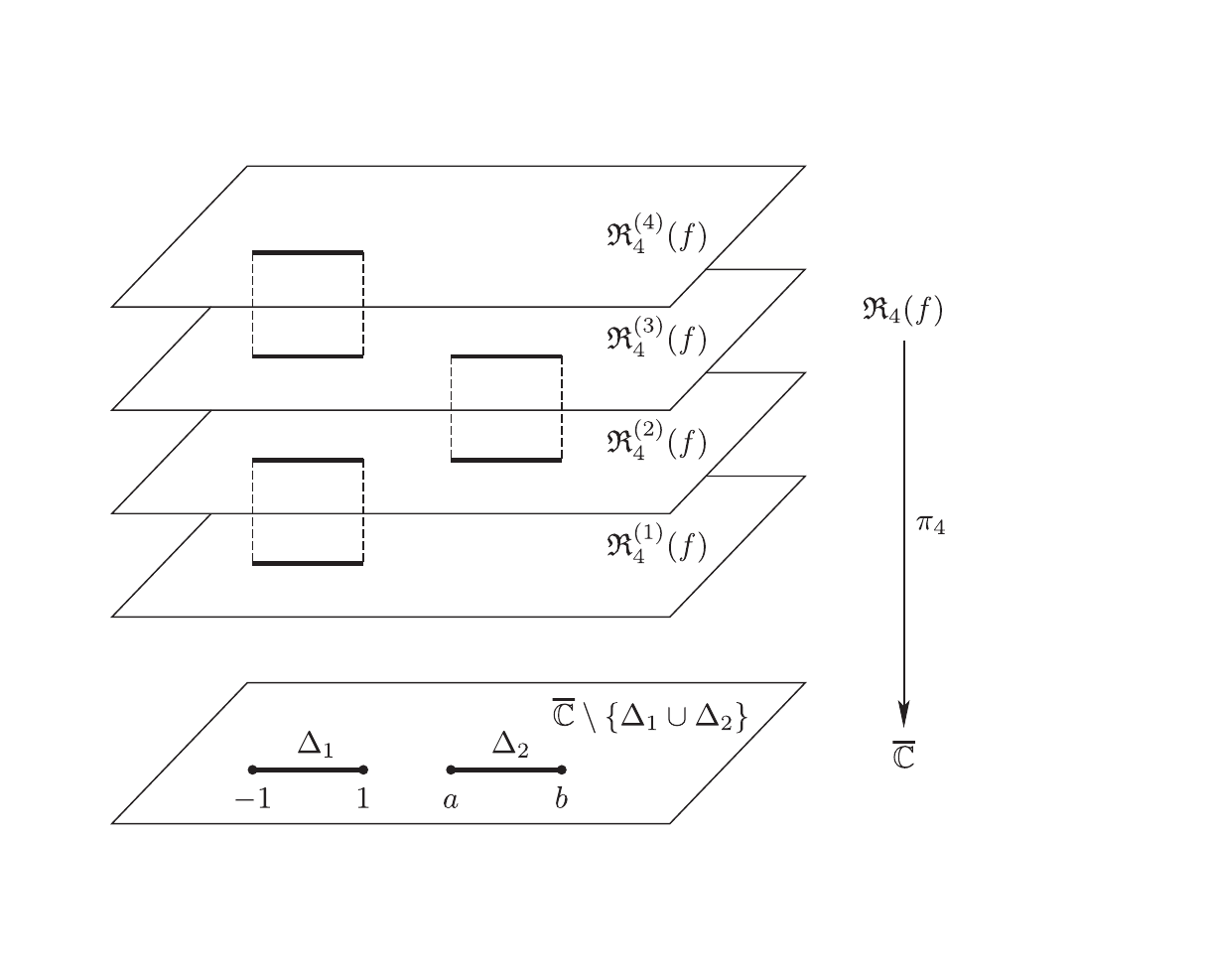}}
\vskip-6mm
\caption{The four-sheeted Riemann surface $\RS_4(\myt{f})$ of the function~$\myt{f}$
given by~\eqref{6}.
}
\label{fig}
\end{figure}

For a~function~$f$ defined by~\eqref{6}, the Riemann surface $\RS(f)$ of the function~$f$
can be looked upon as a~four-sheeted covering of the
Riemann sphere $\myo\CC$ with the branch\footnote{Note that all branch points of this covering are of second order.}
points $z=\pm1,a,b$. Let $\pi_4\colon\RS_4\to\myo\CC$, $\RS_4=\RS_4(f)$
be the corresponding canonical projection. We shall assume that the R.s.\ $\RS_4(f)$
is realized as follows (see Fig.~\ref{fig}). The first (open)
sheet\footnote{Throughout, by
sheets of the R.s.\ we mean open sheets.}
 $\RS^{(1)}_4$ of the R.s.\ $\RS_4(f)$ is the Riemann sphere cutted
along the interval~$\Delta_1$. The second $\RS^{(2)}_4$
and the third $\RS^{(3)}_4$ sheets are the Riemann spheres with cuts, respectively, along the intervals
$\Delta_1$ and $\Delta_2$. Further, the fourth sheet $\RS^{(4)}_4$ is the Riemann sphere with cut
along the interval $\Delta_2$. We follow the standard assumption that each cut along the
corresponding interval has two ``edges'' (the upper and bottom ones).
The four sheets are ``glued'' into the surface by the standard crosswise  identification of the upper
and bottom edges of the corresponding cuts of different sheets.
According to this rule, the first sheet is connected to the second one by the cut corresponding
to the interval $\Delta_1$, after this the second sheet is glued to the third sheet along the
cut corresponding to the interval $\Delta_2$,
and lastly, the third sheet is glued to the fourth sheet along the cut corresponding to the
interval~$\Delta_1$.
The genus of the R.s.\ $\RS_4(f)$ thus obtained is zero, which means that this
R.s.\ is topologically equivalent to the Riemann sphere~$\myo\CC$.

The points lying on the sheets
$\RS_4^{(1)}$, $\RS_4^{(2)}$, $\RS_4^{(3)}$ and $\RS_4^{(4)}$
of this R.s.\  will be denoted by
$z^{(1)},z^{(2)},z^{(3)}$ and $z^{(4)}$, respectively.
We set
$\Gamma^{(1,2)}=\partial\RS_4^{(1)}\cap\partial\RS_4^{(2)}$,
$\Gamma^{(2,3)}=\partial\RS_4^{(2)}\cap\partial\RS_4^{(3)}$,
$\Gamma^{(3,4)}=\partial\RS_4^{(3)}\cap\partial\RS_4^{(4)}$.
For a~function $f\in\SZ$ of the form~\eqref{6}, we have $\pi_4(z^{(j)})=z\in
D_1\setminus\Delta_2$,
$j=1,2,3,4$, $\pi_4(\Gamma^{(1,2)})=\pi_4(\Gamma^{(3,4)})=\Delta_1$,
and $\pi_4(\Gamma^{(2,3)})=\Delta_2$.

In the above notation, the aforementioned results of Stahl's theory
(see \eqref{4}--\eqref{5}) for a~function~$f\in\SZ$ defined by~\eqref{6} can be
interpreted as follows:
\begin{align}
[n/n]_{f_{\infty}}(z)&\to f(z^{(1)}),\quad z\in D_1,
\label{7}\\
\bigl|f(z^{(1)})-[n/n]_{f_\infty}(z)\bigr|^{1/n}
&\to e^{-2g_{\Delta_1}(z,\infty)},\quad z\in D_1
\label{8}
\end{align}
as $n\to\infty$; the convergence in \eqref{7} and~\eqref{8} is uniform inside the domain $D_1$
(i.e., on compact subsets of~$D_1$). Note that in~\eqref{7} and~\eqref{8}
the convergence in capacity can be replaced by the uniform convergence, because
a~function~$f$ of the form~\eqref{6} is easily seen to be a~Markov function; i.e., it can be represented~as
$$
f(z)=\frac1{\sqrt{AB}}+\myh{\sigma}(z),
$$
where
$$
\myh{\sigma}(z):=\int_{-1}^1\frac{d\sigma(x)}{z-x}
$$
is the Cauchy transform of a~measure~$\sigma$ supported on the interval~$\Delta_1$.

Thus, the maximal Stahl domain $D_1$ corresponds to the first sheet
$\RS_4^{(1)}$ of the R.s.\ $\RS_4(f)$, which, as was noted above, is topologically
equivalent to the Riemann sphere~$\myo\CC$. On this Riemann sphere the domain
corresponding to the first sheet occupies, in a~sense, only the fourth part of the sphere.
Therefore, the maximality of the Stahl domain should be understood as the maximality
related only to diagonal Pad\'e approximants: these rational approximants do not converge
on any larger domain.
The question of the efficient analytic continuation of a~given germ~$f_\infty$ to the other
sheets of the R.s.\ $\RS_4(f)$ is shown to be solvable using rational approximants
based on Hermite--Pad\'e polynomials.

Until the end of \S\,\ref{s2} we assume that $f\in\SZ$ is given by
representation~\eqref{6}, $f\in\HH(D_1)$, and $f_\infty\in\HH(\infty)$ is the
corresponding germ defined at the point $z=\infty$.

\subsection{}\label{s2s2}
For the system of three germs $f_\infty,f_\infty^2,f_\infty^3$ and an arbitrary $n\in\NN$, we define the
type~II Hermite--Pad\'e polynomials $Q_{3n},P_{3n,1},P_{3n,2},P_{3n,3}\in\PP_{3n}$,
$Q_{3n}\not\equiv0$, by the relations\footnote{The polynomial $Q_{3n}$ plays the principal role in this construction.
From this polynomial, the remaining polynomials
$P_{3n,j}$, $j=1,2,3$, can be determined uniquely.}
(see \cite{Nut84},~\cite{NiSo88}):
\begin{align}
(Q_{3n}f_\infty-P_{3n,1})(z)&=O\(\frac1{z^{n+1}}\),\quad z\to\infty,
\label{11}\\
(Q_{3n}f_\infty^2-P_{3n,2})(z)&=O\(\frac1{z^{n+1}}\),\quad z\to\infty,
\label{12}\\
(Q_{3n}f_\infty^3-P_{3n,3})(z)&=O\(\frac1{z^{n+1}}\),\quad z\to\infty.
\label{13}
\end{align}
It is easily checked that for the function~$f$ defined in~\eqref{6} the system
$f,f^2,f^3$ is a~Nikishin system (regarding the definition and properties
of Nikishin systems, see \cite{Nik86}, \cite{NiSo88},~\cite{FiLo11}). Namely,
the following representations hold (see~\cite{IkKoSu18}):
\begin{equation}
\begin{gathered}
f(z)=\frac1{\sqrt{AB}}+\myh{\sigma}(z),\quad
f^2(z)=\frac1{AB}+\frac1{\sqrt{AB}}\myh{\sigma}(z)-\myh{s}_1(z),
\\
f^3(z)=\frac1{\sqrt{(AB)^3}}+\frac1{AB}\myh\sigma(z)-\frac1{\sqrt{AB}}
\myh{s}_1(z)-\myh{s}_2(z),\\
\text{where}\quad s_1:=\<\sigma,\sigma_2\>,
\quad s_2:=\<\sigma,\sigma_2,\sigma\>
\end{gathered}
\label{13.2}
\end{equation}
(for the notation used in \eqref{13.2}, see~\cite{GoRaSo97}).
Hence, since a~Nikishin system is perfect (see~\cite{FiLo11}), we have
$\mdeg{Q_{3n}}=3n$, $\mdeg{P_{3n,j}}=3n$, $j=1,2,3$.
Therefore, by normalizing the polynomials $Q_{3n}$ all these four polynomials are determined uniquely.

We introduce now the necessary notation.

By $g_{\Delta_j}(z,\zeta)$ we shall denote the Green function of the domain
$D_j=\myo\CC\setminus\Delta_j$, $j=1,2$, with singularity at the point $z=\zeta$.
The corresponding Green potential of a~measure~$\mu$, $\supp\mu\subset\CC$, is denoted by
$G^\mu_{\Delta_j}(z)$,
$$
G^\mu_{\Delta_j}(z):=\int g_{\Delta_j}(\zeta,z)\,\mu(\zeta),
\quad z\in D_j.
$$
By $M_1(\Delta_j)$ we denote the space of all unit (positive
Borel) measures with support on~$\Delta_j$, $j=1,2$. Further,
$\beta_{\Delta_j}(\mu)$ denotes the balayage of an arbitrary measure~$\mu$
from the domain~$D_j$ onto its boundary $\Delta_j=\partial D_j$.

\begin{lemma}\label{lem1}
In the class $M_1(\Delta_2)$, there exists a~unique measure~$\lambda$ such that
on~$\Delta_2$
\begin{equation}
V^{\lambda}(x_2)+G^{\lambda}_{\Delta_1}(x_2)+g_{\Delta_1}(x_2,\infty)
\equiv\const=\gamma(\lambda),\quad x_2\in\Delta_2,
\label{14}
\end{equation}
with some constant $\gamma=\gamma(\lambda)$.
\end{lemma}

Relation \eqref{14} will be called the \textit{equilibrium relation},
$\lambda$~is the \textit{equilibrium measure}, and $\gamma(\lambda)$ is the \textit{equilibrium  constant}.

\begin{remark}\label{rem1}
It is worth pointing out that the equilibrium relation \eqref{14} holds for the mixed
Green-logarithmic potential $V^\mu(z)+G^\mu_{\Delta_1}(z)$ with the external field
$\psi(z)=g_{\Delta_1}(z,\infty)$;
cf.~\cite{GoRa81}, \cite{GoRa84}, \cite{RaSu13},~\cite{BuSu15}.
\end{remark}

Let $\lambda_1:=\beta_{\Delta_1}(\lambda)\in M_1(\Delta_1)$ be the balayage of the measure~$\lambda$
from~$D_1$ onto~$\Delta_1$,
\begin{equation}
\lambda_2:=\frac13\bigl(2\tau^{\vphantom{p}}_{\Delta_1}+\lambda_1\bigr)
\in M_1(\Delta_1)
\label{14.2}
\end{equation}
(recall that $\tau^{\vphantom{p}}_{\Delta_1}$ is the Robin measure of the interval
$\Delta_1$).

The following result holds.

\begin{theorem}\label{th1}
The degree of the polynomial $Q_{3n}$  is $3n$, all zeros of the polynomial $Q_{3n}$ are simple
and lie in the interval $\Delta_1^\circ:=(-1,1)$.
Moreover, if $n\to\infty$, then
\begin{align}
\frac1n\chi(Q_{3n})&\overset{*}\longrightarrow3\lambda_2;
\label{15}\\
\frac{P_{3n,1}}{Q_{3n}}(z)&\to f(z^{(1)}),\quad z\in D_1;
\label{16}
\end{align}
the convergence rate in \eqref{16} can be characterized by the relation
\begin{equation}
\biggl|f(z^{(1)})-\frac{P_{3n,1}}{Q_{3n}}(z)\biggr|^{1/n}
\to \delta_1(z),\quad z\in D_1,
\label{17}
\end{equation}
where
$$
\delta_1(z)
:=e^{-\bigl(G^{\lambda}_{\Delta_1}(z)+2g_{\Delta_1}(z,\infty)\bigr)}<1,
\quad z\in D_1
$$
{\rm (}cf.~\eqref{7} and~\eqref{8}{\rm )}. The convergence in \eqref{16} and~\eqref{17} is uniform
inside\footnote{That is, on compact subset of the domain;
cf.~\cite{Mar68}.}
the domain~$D_1$.
\end{theorem}

\begin{remark}\label{rem2}
Note that if $n\to\infty$, then
$$
\frac{P_{3n,2}}{Q_{3n}}(z)\to f^2(z^{(1)}),\quad
\frac{P_{3n,3}}{Q_{3n}}(z)\to f^3(z^{(1)}),
\quad z\in D_1.
$$
However, these relations provide no additional information about the values of the function
$f(z^{(1)})$ for $z\in D_1$.
\end{remark}

\subsection{}\label{s2s3}
For the family of germs $[1,f_\infty,f_\infty^2,f_\infty^3]$ and an arbitrary
$n\in\NN$, the type~I Hermite--Pad\'e polynomials $Q_{n,j}\in\PP_n$, $Q_{n,j}\not\equiv0$,
$j=0,1,2,3$, are defined (not uniquely) by the relation\footnote{In this construction
the polynomials $Q_{n,1}$, $Q_{n,2}$, $Q_{n,3}$ play the principal role. From this polynomials
the polynomial $Q_{n,0}$ is determined uniquely.}
\begin{equation}
R_n(z):=(Q_{n,0}+Q_{n,1}f_\infty+Q_{n,2}f_\infty^2+Q_{n,3}f_\infty^3)(z)
=O\(\frac1{z^{3n+3}}\),\quad z\to\infty
\label{18}
\end{equation}
(for more details, see \cite{Nut84},~\cite{NiSo88}).

The following result holds.

\begin{theorem}\label{th2}
For any $n\in\NN$\enskip  $\mdeg{Q_{n,j}}=n$, $j=0,1,2,3$. All zeros of the polynomials
$Q_{n,j}$ are simple and lie in the interval $\Delta_1^\circ$. Moreover, if $n\to\infty$, then
\begin{align}
\frac1n\chi(Q_{n,j})\overset{*}\longrightarrow\lambda_1\in M(\Delta_1),\quad
j=0,1,2,3;
\label{19}\\
\frac{Q_{n,2}}{Q_{n,3}}(z)\to-\bigl[f(z^{(1)})+f(z^{(2)})+f(z^{(3)})\bigr],
\quad z\in D_1;
\label{20}
\end{align}
the convergence rate in~\eqref{20} is characterized by the relation
\begin{equation}
\biggl|
\bigl[f(z^{(1)})+f(z^{(2)})+f(z^{(3)})\bigr]+\frac{Q_{n,2}}{Q_{n,3}}(z)
\biggr|^{1/n}\to\delta_3(z),\quad z\in D_1,
\label{21}
\end{equation}
where
$$
\delta_3(z):=e^{-2G^{\lambda}_{\Delta_1}(z)}<1,\quad z\in D_1
$$
{\rm (}cf.~\eqref{7} and~\eqref{8}{\rm )}. The convergence in \eqref{20} and~\eqref{21}
is uniform inside the domain~$D_1$.
\end{theorem}

\begin{remark}\label{rem3}
Note that convergence in~\eqref{20} follows directly from more general results
of~\cite{ChKoPaSu17}. The description of the measures  $\lambda_2$
and $\lambda_1$ that characterize the limit distribution of the zeros of the polynomials
$Q_{3n}$ and $Q_{n,j}$ in terms of the  {\it scalar} equilibrium problem~\eqref{14}
seems to be new. Usually such a~description is given in terms of the {\it
vector} equilibrium problem (see, \textit{prima  facie}, \cite{GoRa81},~\cite{GoRaSo97},
and also \cite[\S~5]{ApLy10}).
This remark also applies to assertion~\eqref{15} of Theorem~\ref{th1}
on the limit distribution of the zeros of the type~II Hermite--Pad\'e polynomials  $Q_{3n}$.
Namely, as was pointed out above, the system of functions $f,f^2,f^3$ forms a~Nikishin system.
In this case, the vector equilibrium problem with the Nikishin interaction matrix
(see, above all, \cite{Nik86}, and also
\cite{NiSo88},~\cite{GoRaSo97}, \cite[\S~5]{ApLy10})
is traditionally used to describe the limit distribution of the zeros of type~II Hermite--Pad\'e polynomials.
For the system of functions $f,f^2,f^3$ such a~matrix
is a~$3\times3$-matrix. In contrast to the traditional vector approach, we have succeeded in
obtaining a~complete
characterization of  the limit distribution of the zeros of type~I and type~II Hermite--Pad\'e polynomials
in terms of the scalar equilibrium problem~\eqref{14}.
It is quite likely that, for an arbitrary Nikishin system consisting of any
number of functions, the description of the limit distribution of the zeros of the corresponding
Hermite--Pad\'e polynomials can be phrased in terms of an appropriate scalar equilibrium problem.
This assumption pertains, of course, only to the diagonal case. The off-diagonal case
has some special features  (see,
for example,~\cite{LoVa18} and the references given therein).

\end{remark}

\subsection{}\label{s2s4}
In relation \eqref{18} it was assumed that $\mdeg{Q_{n,j}}\leq{n}$ for all $j=0,1,2,3$.
Following the standard approach (see \cite{Nut84},~\cite{NiSo88}), we now introduce
the three-di\-men\-sio\-nal multiindices $\nn_1:=(n,n-1,n-1)$, $\nn_2:=(n,n,n-1)$, and the
corresponding type~I Hermite--Pad\'e polynomials $Q_{\nn_1,j}$ and
$Q_{\nn_2,j}$, $j=0,1,2,3$.

For the multiindex $\nn_1$, the polynomials
$Q_{\nn_1,0},Q_{\nn_1,1}$ lie in $\PP_n^{*}$, the polynomials
$Q_{\nn_1,2},Q_{\nn_1,3}$ lie in~$\PP_{n-1}^{*}$, and moreover
(see~\eqref{14}),
\begin{equation}
R_{\nn_1}(z):=
(Q_{\nn_1,0}+Q_{\nn_1,1}f_\infty+Q_{\nn_1,2}f_\infty^2+Q_{\nn_1,3}f_\infty^3)(z)
=O\(\frac1{z^{3n+1}}\),\quad z\to\infty,
\label{22}
\end{equation}
$\mdeg{Q_{\nn_1,0}}=\mdeg{Q_{\nn_1,1}}=n$,
$\mdeg{Q_{\nn_1,2}}=\mdeg{Q_{\nn_1,3}}=n-1$,
all zeros of the polynomials $Q_{\nn_1,j}$,
$j=0,1,2,3$, are simple and lie in the interval $\Delta_1^\circ$.

For the multiindex $\nn_2$, the polynomials
$Q_{\nn_2,0},Q_{\nn_2,1},Q_{\nn_2,2}$ lie in $\PP_n^{*}$, the polynomial
$Q_{\nn_2,3}$ lies in $\PP_{n-1}^{*}$, and moreover,
\begin{equation}
R_{\nn_2}(z):=
(Q_{\nn_2,0}+Q_{\nn_2,1}f_\infty+Q_{\nn_2,2}f_\infty^2+Q_{\nn_2,3}f_\infty^3)(z)
=O\(\frac1{z^{3n+2}}\),\quad z\to\infty,
\label{23}
\end{equation}
$\mdeg{Q_{\nn_2,0}}=\mdeg{Q_{\nn_2,1}}=\mdeg{Q_{\nn_2,2}}=n$,
$\mdeg{Q_{\nn_2,3}}=n-1$, all zeros of the polynomials $Q_{\nn_2,j}$,
$j=0,1,2,3$, are simple and lie in the interval $\Delta_1^\circ$.

We now set
\begin{equation}
\renewcommand*{\arraystretch}{1.2}
\begin{aligned}
S_{2n,1}(z):&=
\begin{vmatrix}
Q_{\nn_1,1}(z)&Q_{\nn_1,3}(z)\\Q_{\nn_2,1}(z)&Q_{\nn_2,3}(z)
\end{vmatrix}
=\bigl(Q_{\nn_1,1}Q_{\nn_2,3}-Q_{\nn_1,3}Q_{\nn_2,1}\bigr)(z)\in\PP_{2n-1},\\
S_{2n,2}(z):&=
\begin{vmatrix}
Q_{\nn_1,2}(z)&Q_{\nn_1,3}(z)\\Q_{\nn_2,2}(z)&Q_{\nn_2,3}(z)
\end{vmatrix}
=\bigl(Q_{\nn_1,2}Q_{\nn_2,3}-Q_{\nn_1,3}Q_{\nn_2,2}\bigr)(z)\in\PP_{2n-1}.
\end{aligned}
\label{24}
\end{equation}

The following result holds.

\begin{theorem}\label{th3}
The interval $\Delta_2$ attracts as $n\to \infty$ all zeros of the polynomials $S_{2n,j}$, $j=1,2$,
and moreover,
\begin{align}
\frac1{2n}\chi(S_{2n,j})&\overset{*}\longrightarrow\lambda\in M_1(\Delta_2);
\label{25}\\
\frac{S_{2n,1}}{S_{2n,2}}(z)&\overset\mcap\longrightarrow -\bigl[f(z^{(1)})+f(z^{(2)})\bigr],
\quad z\in D_2=\myo\CC\setminus\Delta_2;
\label{26}
\end{align}
the convergence rate in \eqref{26} can be characterized by the relation
\begin{equation}
\biggl|\bigl[f(z^{(1)})+f(z^{(2)})\bigr]+\frac{S_{2n,1}}{S_{2n,2}}(z)\biggr|^{1/n}
\overset\mcap\longrightarrow \delta_2(z),\quad z\in D_2,
\label{27}
\end{equation}
where
$$
\delta_2(z):=e^{2\(V^{\lambda}(z)+G^{\lambda}_{\Delta_1}(z)+g_{\Delta_1}(z,\infty)
-\gamma(\lambda)\)}<1,\quad z\in D_2.
$$
\end{theorem}

Note that in \eqref{26}--\eqref{27} we assert the convergence in capacity, rather than
the uniform convergence as in \eqref{16}--\eqref{17} and \eqref{20}--\eqref{21}.
However, this constraint,
which depends only on our method of the proof of Theorem~\ref{th3}, is not related to the essence of the matter,
The crux here is that, although
various extensions of Hermite--Pad\'e polynomials are well-known
(see \cite{Sor83},~\cite{FiLoLoSo10}), it seems that
the polynomials $S_{2n,j}$ were not considered earlier and their properties have not been not studied.
Relations \eqref{24} are explicit representations of the polynomials
$S_{2n,j}$ in terms of type~I Hermite--Pad\'e polynomials corresponding to two neighbouring
multiindices. There is a~natural formal definition of these polynomials which is
quite similar to the definitions~\eqref{11}--\eqref{13} and~\eqref{18}
of type~I and~II Hermite--Pad\'e polynomials, respectively.
Such definition will be given in the paper \cite{IkKoSu18}.
From this definition and the above Theorem~\ref{th3} it follows that the
polynomials $S_{2n,j}$ are intermediate between
the type~I Hermite--Pad\'e polynomials and the type~II Hermite--Pad\'e polynomials.

For a~meaningful definition of the polynomials $S_{2n,1}$ and $S_{2n,2}$ in terms of
determinants~\eqref{24}, whose entries are type~I Hermite--Pad\'e polynomials,
it is necessary that these determinants be not identically zero.
To satisfy this requirement, the multiindices $\nn_1$ and $\nn_2$ should obey a~condition
more restrictive than the normality condition for the standard PA. In case of~\eqref{6} (i.e., for real $A$
and~$B$), this condition is satisfied. However, in a~more general setting with complex $A$ and~$B$,
this condition ceases to be true (see \S~\ref{s3} below).

In this connection, we recall the following well-known fact
(see \cite{Nut84},~\cite{NiSo88}).

Consider the multiindices
$\nn_1'=(n,n-1,n-1)$, $\nn_2'=(n-1,n,n-1)$ and $\nn_3'=(n-1,n-1,n)$ (the multiindex
$\nn_1'$ coincides with the above multiindex~$\nn_1$).
Let $Q_{\nn_k',j}$, $j=0,1,2,3$, $k=1,2,3$, be the corresponding type~I Hermite--\allowbreak Pad\'e polynomials.
Then we have the following explicit representations for the type~II
Hermite--\allowbreak Pad\'e polynomials defined by relations \eqref{11}--\eqref{13}
\begin{equation}
\renewcommand*{\arraystretch}{1.2}
\begin{aligned}
Q_{3n}(z)&=
\begin{vmatrix}
Q_{\nn_1',1}(z)&Q_{\nn_1',2}(z)&Q_{\nn_1',3}(z)\\
Q_{\nn_2',1}(z)&Q_{\nn_2',2}(z)&Q_{\nn_2',3}(z)\\
Q_{\nn_3',1}(z)&Q_{\nn_3',2}(z)&Q_{\nn_3',3}(z)
\end{vmatrix},\\
P_{3n,1}(z)&=-
\begin{vmatrix}
Q_{\nn_1',0}(z)&Q_{\nn_1',2}(z)&Q_{\nn_1',3}(z)\\
Q_{\nn_2',0}(z)&Q_{\nn_2',2}(z)&Q_{\nn_2',3}(z)\\
Q_{\nn_3',0}(z)&Q_{\nn_3',2}(z)&Q_{\nn_3',3}(z)
\end{vmatrix}
\end{aligned}
\label{28}
\end{equation}
under the assumption that the determinants do not vanish identically.
Similar representations under the same nondegeneracy conditions also hold for the
polynomials $P_{3n,2}$ and $P_{3n,3}$.

We emphasize once again that definition~\eqref{11}--\eqref{13} of a~type~II Hermite--Pad\'e polynomials
is much more general than the explicit representation~\eqref{28}.
Namely, according to \eqref{11}--\eqref{13} such polynomials $Q_{3n}\not\equiv0$,
$P_{3n,j}\not\equiv0$ always exist, whereas representation \eqref{28} does not always hold.
The situation with the above polynomials
$S_{2n,j}$ is quite similar. As was already mentioned above,
we plan to discuss these polynomials in detail in the separate paper \cite{IkKoSu18}.
Nevertheless, we mention however that from the above it follows that
the type~I Hermite--Pad\'e polynomials themselves play in a~sense a~key role in the efficient solution
of the problem of analytic continuation of a~given germ of a~multivalued
analytic function. Indeed, all rational approximants involved in this procedure
can in principle be calculated in terms of such polynomials corresponding to various multiindices.

\subsection{}\label{s2s3}
Thus, relations \eqref{16},~\eqref{26} and~\eqref{20} collectively compose the
set from which one can recover in succession, from local data
(a~germ~$f_\infty$),
the values of a~multivalued function
$f(z)$ on the first, second, and third sheets of the R.s.~$\RS_4(f)$.
Of course, here we speak about the limit procedure,  in which, for a~given germ $f_\infty\in\HH(\infty)$,
the number of Laurent coefficients involved in the construction of appropriate rational approximants
grows unboundedly, because an analytic function cannot be given by a~finite number of its
Laurent coefficients.
Note that for each $n\in\NN$ the number of Laurent coefficients
involved in the construction of the three functions under consideration is at most $3n+4$.

It is also worth pointing out that $f(z)$ cannot be extended to the last fourth sheet $\RS^{(4)}_4$
of the R.s.\ $\RS_4(f)$  by the above procedure. This is in full accord
with the fact that the diagonal Pad\'e approximants extend the germ $w_\infty\in\HH(\infty)$
of the hyperelliptic function~$w$,
$w^2=P_{2\ell}(z)=z^{2\ell}+\dotsb\in\PP_{2\ell}$, only on the first sheet
of the corresponding two-sheeted hyperelliptic surface. This seems to be quite natural,
because, in simple terms, the values of the hyperelliptic function~$w$ differ on two sheets of the
corresponding R.s.\ only by the sign.

We also emphasize the following fact.
The reader may well have formed the impression that
on the first sheet
$\RS^{(1)}_4$ of the R.s.\ $\RS_4(f)$ the function~$f$ can be always recovered from the germ  $f_\infty$
using the diagonal PA  $[n/n]_{f_\infty}$. But this is not so. The thing is that when speaking about the
partition of the R.s.\ $\RS_4(f)$ into ``sheets'' we  have tacitly meant an intuitive
(but never formalized) perception about the sheets of the R.s.\ $\RS_4(f)$
based on the cuts drawn along the already specified intervals
$\Delta_1=[-1,1]$ and $\Delta_2=[a,b]$. But in the general case, even for a~simple function of the form~\eqref{6},
this is no longer the case and the thing is more complicated.
Namely,
the situation changes drastically if at least one of the two points $a=(A+1/A)/2$ or $b=(B+1/B)/2$
moves from the real line into the complex plane:
the natural (for the Hermite--Pad\'e polynomials) partition of the R.s.\ $\RS_4(f)$ into sheets will proceed
not along the original intervals $\Delta_1$ and $\Delta_2$, but along some arcs
$\ell_1$ and $\ell_2$ that joint the points $-1,1$ and $a,b$, respectively.
More precisely, here we speak about the so-called {\it canonical} or the {\it Nuttall} partition of the R.s.\ into sheets,
which is uniquely determined by the real part of some Abelian integral with logarithmic
singularities and purely imaginary periods (see
\cite[Sect.~3]{Nut84}, \cite[formula (3) and Lemma~5]{ChKoPaSu17}, and also \S\,\ref{s3}
below). In particular, for such partition into sheets, to the first sheet there corresponds a~domain
$\myo\CC\setminus\ell_1$, which in general is distinct from the domain
$D_1=\myo\CC\setminus\Delta_1$ corresponding to the first sheet of the Stahl surface.

It can be shown that in the case when the quantities $A$ and $B$ in~\eqref{6} are real,
the Nuttall partition of the  R.s.\ $\RS_4(f)$ into sheets
is in accord with the intuitive perception about this partition.

Indeed, it can be easily shown that the real function~$u$,
\begin{equation}
\begin{aligned}
u(z^{(1)}):&=3V^{\lambda_2}(z)\quad z\in D_1,\\
u(z^{(2)}):&=2V^{\lambda}(z)-3V^{\lambda_2}(z)+2c_1,\quad
z\in\myo\CC\setminus(\Delta_1\cup\Delta_2 ),\\
u(z^{(3)}):&=V^{\lambda_1}(z)-2V^{\lambda}(z)+2c_1+c_2,\quad
z\in\myo\CC\setminus(\Delta_1\cup\Delta_2 ),\\
u(z^{(4)}):&=-V^{\lambda_1}(z)+2c_1+c_2+2c_3,\quad z\in D_1,
\end{aligned}
\label{29}
\end{equation}
which is defined in terms of the above (see \eqref{14} and~\eqref{14.2}) measures $\lambda$, $\lambda_1$ and
$\lambda_2$,
where $c_1,c_2,c_3$ are appropriately  chosen constants, is a~(singlevalued) harmonic function
on the R.s.\ $\RS_4(f)$ everywhere except the four points $\pi_4^{-1}(\infty)$
at which it has logarithmic singularities:
\begin{equation}
\begin{aligned}
u(z^{(1)})&=-3\log|z|+O(1),\quad z\to\infty,\\
u(z^{(j)})&=\log|z|+O(1),\quad z\to\infty,\quad j=2,3,4.
\end{aligned}
\notag
\end{equation}
Moreover,
\begin{equation}
u(z^{(1)})<u(z^{(2)})<u(z^{(3)})<u(z^{(4)}),
\quad z\in\myo\CC\setminus(\Delta_1\cup\Delta_2).
\label{30}
\end{equation}
But inequalities~\eqref{30} mean that  our adopted (on an intuitive level) partition of
the  R.s.\ $\RS_4(f)$ into sheets is in fact the Nuttall partition
\cite[Sect.~3, relation (3.1.3)]{Nut84}; see
also~\cite{Apt99},~\cite{MaRaSu16} and Fig.~\ref{fig}.

The next section will be dedicated to the complex case\,--\,in this setting
the dominant role is now played the Nuttall partition of the  R.s.\ $\RS_4(f)$ into sheets.
The properties of this partition in the complex case are far from being evident
\cite[Sect.~3]{Nut84},
\cite[Lemma~5]{ChKoPaSu17}.

\section{The complex case}\label{s3}

\subsection{}\label{s3s1}
In this section, we formulate and discuss some theoretical results which we have obtained in the
case when in~\eqref{1} all  exponents
$\alpha_j$ are equal to $\pm1/2$ and $m$~is an even number. The complex quantities
$A_j$ are subject to no other constraints than the condition that they are pairwise distinct
and $|A_j|>1$ for any~$j$.
In this setting, for such functions~$f$ from the class~$\SZ$ with the four-sheeted R.s.,  we shall
establish some properties of the Nuttall partition of this R.s.\ $\RS_4(f)$;
such properties do not hold in a~general case\footnote{The general Nuttall conjecture~\cite[Sect.~3]{Nut84}
stating that under the canonical partition into sheets of any $m$-sheeted surface the
complement of the last $m$th sheet is always connected was recently proved in \cite[Lemma~5]{ChKoPaSu17}.}.
Moreover, in this particular case, our alternative approach
to the characterization of the Nuttall partition of the R.s.\ $\RS_4(f)$ into sheets
is capable of delivering some additional information about the properties of such a~partition.

In the case $m=2$ (under the same conditions on the exponents $\alpha_{1,2}=\pm1/2$)
we shall formulate analogues of Theorems~\ref{th1}--\ref{th3}, which hold, in contrast to the real case,
under some additional assumptions on the remainder functions $R_{\nn_1}$, $R_{\nn_2}$, and~$R_{\nn_3}$.
It should be noted that these constrains on the remainder functions
are not related to the essence of the matter, but are rather due
to the imperfection of the currently available methods of dealing with
asymptotical  properties of the Hermite--Pad\'e polynomials. In particular, many
conjectures in this direction (see \cite{Nut84}, \cite{Sta88}, \cite{Apt08}) remain unjustified at present.

First fairly general results on the convergence (in capacity)
of diagonal PA were obtained by Nuttall (see \cite{NuSi77},~\cite{Nut84})
in the class of hyperelliptic functions (in particular, these are the functions whose branch points
are only of the second order). It is in this case in which in the early 1980s (i.e., before Stahl's works of 1985--1986)
the Nuttall partition of the hyperelliptic surface $\RS_2(w)$ into sheets was
introduced using the real part of an Abelian integral
of the third kind with purely imaginary periods and appropriate logarithmic singularities
at the two points $\pi^{-1}_2(\infty)$
(see also~\eqref{32} and~\eqref{33}). In this particular case, each of the two Nuttall sheets is a~domain.
In 1984 Nuttall~\cite[Sect.~3]{Nut84} put forward the conjecture that for any $m\ge3$
and for the partition of an $m$-sheeted R.s.\ into Nuttall sheets, the complement of the closure of the
last ``uppermost'' (cf.~\eqref{30}) sheet is a~domain on the original R.s.\ $\RS_m$.
In its complete form this conjecture was proved very recently (see Lemma~5 of~\cite{ChKoPaSu17}).
In the present paper, for the above class of multivalued analytic functions
of the form~\eqref{1}, we prove, under the condition $\alpha_j=\pm1/2$, some additional properties of
such a~partition; in particular, we show that the first sheet $\RS^{(1)}_4(f)$ is connected.

\subsection{}\label{s3s2}
So, let
\begin{equation}
f(z):=\prod_{j=1}^m\(A_j-\frac1{\pfi(z)}\)^{\alpha_j},\quad z\in D_1,
\label{31}
\end{equation}
where $\alpha_j=\pm1/2$ for any~$j$, $m=0\pmod{2}$, all $A_j$ are pairwise distinct, and $|A_j|>1$.
The class of such functions will be denoted by $\SZ_{1/2}$.

Consider the two-sheeted R.s.\ $\RS_2(w)$ of the function~$w$ defined by the equality
$w^2=z^2-1$.
We shall assume that the R.s.\ $\RS_2(w)$ is realized as a~two-sheeted covering
of the Riemann sphere $\myo\CC$ using the explicitly given uniformization
\begin{equation}
z=\frac12\(\zeta+\frac1\zeta\),
\quad
w=\frac12\(\zeta-\frac1\zeta\),\quad \zeta\in\myo\CC_\zeta.
\label{32.0}
\end{equation}
Accordingly, the first sheet $\RS^{(1)}_2$, on which
$w=(z^2-1)^{1/2}/z\to1$ as $z\to\infty$, corresponds to the exterior $\{\zeta:|\zeta|>1\}$ of the unit disc
$\DD_\zeta:=\{\zeta:|\zeta|<1\}$ in the $\zeta$-plane, and the second sheet
$\RS^{(2)}_2$, on which $w=-(z^2-1)^{1/2}/z\to-1$ as $z\to\infty$, corresponds to the unit disc itself.
By a~point~$\zz$ on the R.s.\ $\RS_2(w)$, $\zz\in\RS_2(w)$, we shall mean the pair
$\zz:=(z,w)=(z,\pm(z^2-1)^{1/2})$. The canonical projection
$\pi_2\colon\RS_2(w)\to\myo\CC$ is defined by $\pi_2(\zz):=z$.
By the point $\zz=\infty^{(1)}\in\RS_2(w)$ we mean the point on the R.s.\ $\RS_2(w)$
such that $\pi_2(\infty^{(1)})=\infty$ and $w/z\to1$ as $\zz\to\infty^{(1)}$. Similarly, for $\zz=\infty^{(2)}$ we have
$\pi(\infty^{(2)})=\infty$ and $w/z\to-1$ as $\zz\to\infty^{(2)}$. A~passage  from the first sheet $\RS^{(1)}_2$ to the
second sheet $\RS^{(2)}_2$ proceeds along the cut interval $\Delta_1$ (here we assume as usual that
an interval has two edges, the upper and bottom ones) and that the sheets are  ``glued'' by
identifying crosswisely the edges of cuts; i.e.,
by identifying the upper edge of one cut with the lower edge of the other, and vice versa.
It can be easily shown that the above partition into sheets
is a~Nuttall partition. Indeed, let
\begin{equation}
G(\zz):=\int_{-1}^{\zz}\frac{dt}{\sqrt{t^2-1}}=\log(z+w)=\log(z\pm(z^2-1)^{1/2})
\label{32}
\end{equation}
be an Abelian integral of the third kind with purely imaginary periods and logarithmic singularities
only at the points  $\zz=\infty{(1)},\infty^{(2)}$. Then
$u_2(\zz):=-\Re G(\zz)$ is a~harmonic function on the R.s.\ $\RS_2(w)\setminus\{\infty^{(1)}, \infty^{(2)}\}$,
$$
u_2(\zz)=\mp\log|z|+O(1),\quad \zz\to\infty^{(1)},{\infty^{(2)}},
$$
and
\begin{equation}
u_2(z^{(1)})<u_2(z^{(2)}),\quad z\in D_1.
\label{33}
\end{equation}
Moreover, $u_2(\zz)=0$ for $\zz\in\bGamma$, $\bGamma:=\pi^{-1}_2(\Delta_1)$.
So,
$-u_2(z^{(1)})=g_{\Delta_1}(z,\infty)$ is the Green function of the domain~$D_1$.
An arbitrary germ $f_\infty\in\HH(\infty)$ of the function $f\in\SZ_{1/2}$,
$f\in\HH(D_1)$, is lifted to the point $\zz=\infty^{(1)}\in\RS_2(w)$ and extends
to the entire first sheet $\RS^{(1)}_2(w)$ as a~singlevalued  holomorphic function
(recall that $\pi_2(\RS_2^{(1)}(w))=D_1=\myo\CC\setminus\Delta_1$). A~further
singlevalued extension of this function to the entire second sheet of the R.s.\ $\RS_2(w)$
is hindered by the branch points $a^{(2)}_j\in\RS^{(2)}_2(w)$ such that
$\pi_2(a^{(2)}_j)=a_j$, $j=1,\dots,m$. In order that such a~singlevalued
analytic (meromorphic) extension of the germ $f_\infty$ be possible one should make appropriate cuts
on the second sheet $\RS_2^{(2)}(w)$.

\begin{definition}[\rm (cf.~\cite{Sue17}, Definition~2)]\label{def1}
A~compact set $\KK=K^{(2)}\subset\RS_2^{(2)}(w)$ is called {\it admissible} for a~germ
$f_\infty\in\HH(\infty^{(1)})$ of a~function $f\in\SZ_{1/2}$ if this germ
$f_\infty$ extends as a~singlevalued meromorphic  function to the domain
$\Omega_{\infty^{(1)}}(K^{(2)})$, which is the connected component of the complement of~$K^{(2)}$
containing the point $\zz=\infty^{(1)}$,
$f_\infty\in\MM(\Omega_{\infty^{(1)}}(K^{(2)}))$.
\end{definition}

The class of all admissible compact sets for a~germ $f_\infty$ will be denoted by
$\myk(f_\infty)$.

Note that one can slightly reduce the class of admissible compact sets and
consider \textit{a~priori} only compact sets with connected complement; i.e.,
such that $\Omega_{\infty^{(1)}}(K^{(2)})=\RS_2(w)\setminus K^{(2)}$.
This change has not effect on subsequent considerations. In what follows, we shall assume that the
family $\myk(f_\infty)$ consists of only regular compact sets in the sense that the domain
$\Omega_{\infty^{(1)}}(K^{(2)})$ is regular with respect to the solution of the Dirichlet problem
(see \cite{ScSp54},~\cite{Chi06}).

We fix a~local coordinate $\zeta$ at the point $\zz=\infty^{(1)}$,
$\zeta(\infty^{(1)})=0$. For an arbitrary compact set $\KK\in\myk(f_\infty)$,
consider the Green function $g_{\KK}(\zz,\infty^{(1)})$ of the domain
$\Omega_{(1)}(\KK)$ with logarithmic singularity at the point $\zz=\infty^{(1)}$
with due account of the choice of the local coordinate at this point; i.e., by definition
$g_{\KK}(\zz,\infty^{(1)})=0$ for $\zz\in\partial\Omega_{\infty^{(1)}}(\KK)$
and
\begin{equation}
g_{\KK}(\zz,\infty^{(1)})=\log\frac1{|\zeta|}+\gamma+o(1),
\quad \zz\to\infty^{(1)},
\label{34}
\end{equation}
where $\gamma=\gamma(\KK)$ is the Robin constant for the domain
$\Omega_{\infty(1)}(\KK)$ at the point $\zz=\infty^{(1)}$ with respect to the chosen
local coordinate~$\zeta$ (see \cite[Ch.~8, \S~4]{ScSp54}).

The following lemma follows from general results of~\cite{ScSp54}.

\begin{lemma}\label{lem2}
The class $\myk(f_\infty)$ contains a~unique compact set $\FF$ such that
\begin{equation}
\gamma(\FF)=\max_{\KK\in\myk(f_\infty)}\gamma(\KK).
\label{35}
\end{equation}
This {\it extremal} compact set $\FF$ does not split the R.s.\ $\RS_2(w)$ (i.e.,
$\Omega^{(1)}(\FF)=\RS_2(w)\setminus\FF$ is a~domain), consists of a~finite number of
analytic arcs, and has the following $S$-property:
\begin{equation}
\frac{\partial g_{\FF}(\zz,\infty^{(1)})}{\partial n^{+}}=
\frac{\partial g_{\FF}(\zz,\infty^{(1)})}{\partial n^{-}},
\quad \zz\in\FF^\circ;
\label{36}
\end{equation}
here $\FF^\circ$ is the union of open arcs whose closures constitute~$\FF$,
$\partial/{\partial n^{+}}$ and $\partial/{\partial n^{-}}$ are the
normal derivatives to~$\FF^\circ$  at the point~$\zz$ taken in opposite directions to the sides of~$\FF$
(the validity of~\eqref{36} does not depend on the local coordinate chosen at the point $\zz\in\FF^\circ$).
\end{lemma}

Lemma~\ref{lem2} is proved using  Schiffer's method of interior variations
in accordance with the general idea proposed in \cite[Ch.~8, \S~4]{ScSp54}.

The quantity  $e^{-\gamma(\KK)}$ is naturally called the capacity of a~compact set~$\KK$
(this quantity depends on the local coordinate chosen at the point $\zz=\infty^{(1)}$). So,
$\FF\in\myk(f_\infty)$ is an admissible compact set of minimal capacity (this property is already independent of the choice of a~local coordinate).

In the actual fact, the genus of the R.s.\ $\RS_2(w)$ is zero, and hence Lemma~\ref{lem2}
can be reduced to the planar case and the corresponding compact set of minimal capacity
on the plane. Indeed, a~uniformization of $\RS_2(w)$ is defined
using the Zhukovskii function as follows:
\begin{equation}
z=\frac12\(\zeta+\frac1\zeta\),\quad
w=\frac12\(\zeta-\frac1\zeta\),\quad \zeta\in\myo\CC_\zeta.
\label{37}
\end{equation}
To the first sheet of  $\RS^{(1)}_2(w)$ there corresponds the exterior
$\myo\CC\setminus\myo\DD$ of the unit disc $\DD:=\{\zeta:|\zeta|<1\}$, while the
second sheet is the unit disc~$\DD$. The function $f(z)\in\SZ_{1/2}$
of the form~\eqref{31} is transformed to the function $\myt{f}(\zeta)$ of the form
\begin{equation}
\myt{f}(\zeta)=\prod_{j=1}^m\(A_j-\frac1{z+w}\)^{\alpha_j}
=\prod_{j=1}^m\(A_j-\frac1\zeta\)^{\alpha_j},
\label{38}
\end{equation}
where all $\alpha_j=\pm1/2$. All
$A_j\in\myo\CC\setminus\myo\DD$, and hence all singular points of the function~$\myt{f}$ are of the form
$\zeta=\myt{a}_j=1/A_j\in\DD$. So, to the extremal compact set~$\FF$ for the function
$f$ on the R.s.\ $\RS_2(w)$ there corresponds an admissible compact set of minimal capacity $\myt{F}$ for the function
$\myt{f}$, as given by the element $\myt{f}_\infty\in\HH(\infty)$ (note that
$\myt{f}\in\MM(\myo\CC\setminus\myt{F})$ and at the point $\zeta=0$ the function $\myt{f}$
has a~pole). By well-known properties of a~compact set of minimal capacity,
$\myt{F}\subset\DD$ (or, more precisely, the compact set $\myt{F}$ lies in the convex hull
of the set  $\{\myt{a}_j=1/A_j$, $j=1,\dots,m\}$), $\myt{F}$~consists of a~finite number of
analytic arcs (trajectories of quadratic differentials), $\myt{F}$~does not split the complex plane and
has the $S$-property~\eqref{36}.
Moreover,
\begin{equation}
\myt{F}=\biggl\{\zeta\in\CC:\Re\int_{\myt{a}_j}^\zeta\sqrt{\frac{V_{m-2}(t)}
{A_m(t)}}\,dt=0\biggr\},
\label{39}
\end{equation}
where $A_m(t):=\prod_{j=1}^m(t-\myt{a}_j)$, $V_{m-2}(t):=(t-v_1)\dots(t-v_{m-2})$ is the  corresponding Chebotarev
polynomial, $v_j$, $j=1,\dots,m$, are points of the Chebotarev compact set~$\myt{F}$.
All these properties of a~compact set of minimal capacity are well known
in a~general case and were obtained by Stahl already in~1985
(see \cite{Sta97b}, \cite{Sta12}, as well as \cite{ApBuMaSu11} and~\cite{Sue15}).

In the case under consideration here, all branch points of the function $\myt{f}$
are of the second order. In this setting, the existence and description of
a~compact set of minimal capacity were proved already by Nuttall  \cite{NuSi77},~\cite{Nut84}. In particular,
in the case of general position, all zeros of the polynomial $V_{m-2}$ have even multiplicity, all $v_j\neq
a_k$ for $j=1,\dots,m-2$, $k=1,\dots,m$, and the compact set
$\myt{F}$ consists of $m/2\in\NN$ disjoint analytic arcs that joint pairwisely the points of the set
$\{\myt{a_j},j=1,\dots,m\}$.

Thus, on the second sheet of the R.s.\ $\RS_2(w)$ we have obtained a~system of analytic arcs
which compose the compact set~$\FF$, do not split the R.s.\ $\RS_2(w)$, have the
$S$-property~\eqref{36}, and are such that $f_\infty\in\MM(\RS_2(w)\setminus{\FF})$.

We now proceed as follows. Given $\zz\in\RS_2(w)\setminus\FF$, consider the function
$$
u(\zz):=-2g_{\FF}(\zz,\infty^{(1)})-g(\zz),
$$
where $g(\zz)=\log|z+w|=\log|\zeta|$. The function $v(\zz)$ is harmonic in the domain
$\RS_2(w)\setminus\FF$ except for the points at infinity
$\zz=\infty^{(1)}$ and $\zz=\infty^{(2)}$, where it behaves as
\begin{equation}
u(\zz)=
\begin{cases}
-3\log|z|+O(1),&\zz\to\infty^{(1)},\\
\log|z|+O(1),&\zz\to\infty^{(2)}
\end{cases}
\label{40}
\end{equation}
(here and in what follows we assume for simplicity that $\infty^{(2)}\not\in\FF$;
otherwise a~more accurate analysis is required quite similar to that of~\cite{ChKoPaSu17}).

We set
\begin{equation}
\EE:=\bigl\{\zz\in\RS_2(w):u(z^{(1)})=u(z^{(2)})\bigr\};
\label{41}
\end{equation}
here by $z^{(1)}\in\RS^{(1)}_2(w)$ and $z^{(2)}\in\RS^{(2)}_2(w)$ we mean, as before,
the Nuttall partition of the R.s.\ $\RS_2(w)$ into sheets.
It is easily checked that  $\EE$~is a~closed arc on the R.s.\ $\RS_2(w)$ passing
through the points $\zz=\pm1$, not crossing the compact set $\FF$, and splitting $\RS_2(w)$
into two domains, of which one contains~$\FF$, and the other one, the point
$\zz=\infty^{(1)}$. Let us now take the second copy of the R.s.\ $\RS_2(w)$ (which we denote by $\RS_2'(w)$)
with cuts drawn along the arcs of the compact set~$\FF$, and glue it to the first
copy with the same cuts using the standard crosswise identification of the opposite edges of a~cut.
This gives us a~four-sheeted R.s.\ on which the original germ
$f_\infty\in\HH(\infty^{(1)})$ extends as a~singlevalued meromorphic
function, because this four-sheeted R.s.\ coincides with the R.s.\ $\RS_4(f)$. We have
$g_\FF(\zz,\infty^{(1)})=0$ for $\zz\in\FF$, and hence the function
$g_\FF(\zz,\infty^{(1)})$
extends continuously to the glued R.s.\ $\RS_2'(w)$ by continuity with the sign change.
The $S$-property~\eqref{36} guarantees that this extension is a~harmonic function on the  four-sheeted R.s.\ thus obtained.
The function $g(\zz)$ is harmonic in a~neighbourhood of the compact set~$\FF$, and hence $f$~extends to
the second copy of the R.s.\ $\RS_2(w)$ by duplication of its values from the first copy
of~$\RS_2(w)$. By $\EE'$ we denote the closed curve lying on the second copy of~$\RS_2(w)$ and
corresponding to the curve~$\EE$; by~$\infty^{(4)}$ we denote the point lying on the second copy of~$\RS_2(w)$ and
corresponding to the point $\zz=\infty^{(1)}$.

We now partition the four-sheeted R.s.\ $\RS_4$ thus obtained into
sheets as follows. The first sheet $\RS^{(1)}_4$ is the domain containing the
point $\zz=\infty^{(1)}$ with boundary $\partial\RS^{(1)}_4=\EE$. The second sheet
$\RS^{(2)}_4$ is the domain with boundary $\EE\cup\FF$. The third sheet $\RS^{(3)}_4$
is the domain with boundary $\FF\cup\EE'$. And the fourth sheet $\RS^{(4)}_4$
is the domain with boundary $\EE'$. Given $\zz\in\RS_2'(w)$, we now set
\begin{equation}
u(\zz):=-2g_{\FF}(\zz,\infty^{(1)})-g(\zz);
\label{42}
\end{equation}
the meaning of the functions $g_{\FF}(\zz,\infty^{(1)})$
and $g(\zz)$ for $\zz\in\RS_2'(w)$ was already explained above. It is easily checked that  $u(\zz)$~is
a~harmonic function in the domain
$\RS_4(f)\setminus\{\infty^{(1)},\infty^{(2)},\infty^{(3)},\infty^{(4)}\}$.
Moreover,
\begin{equation}
u(\zz)=
\begin{cases}
-3\log|z|+O(1),&\zz\to\infty^{(1)},\\
\log|z|+O(1),&\zz\to\infty^{(2)},\\
\log|z|+O(1),&\zz\to\infty^{(3)},\\
\log|z|+O(1),&\zz\to\infty^{(4)},
\end{cases}
\label{42.2}
\end{equation}
and besides, for the above partition of the R.s.\ $\RS_4(f)$ into sheets,
\begin{equation}
u(z^{(1)})<u(z^{(2)})<u(z^{(3)})<u(z^{(4)}).
\label{43}
\end{equation}
So, the partition of the R.s.\ $\RS_4(f)$ into sheets, which we introduced using the Green
function $g_{\FF}(\zz,\infty)$ corresponding to the compact set of minimal capacity
$\FF\subset\RS_2(w)$, turned out to be a~Nuttall partition; for more details, see~\cite{IkKoSu18}.

\subsection{}\label{s3s3}
In the remaining part of the present section we shall content ourselves with the case $m=2$
in representation~\eqref{31}. To be more precise, we shall assume as in \S~\ref{s2}  that
$f\in\SZ$ has the form
\begin{equation}
f(z)=\biggl[\(A-\frac1{\pfi(z)}\)\(B-\frac1{\pfi(z)}\)\biggr]^{-1/2},
\label{44}
\end{equation}
where, as before, $|A|>1$, $|B|>1$, but in contrast to~\S~\ref{s2},
$A$ and~$B$ are not supposed to be real.

We fix a~germ $f_\infty\in\HH(\infty)$ of a~function $f\in\HH(D_1)$
of the form \eqref{44} (two possible germs differ only by the sign). By the
original assumption about the inverse of the Zhukovskii function (namely,
$|\pfi(z)|>1$ for $z\notin\Delta_1$) we have $f_\infty\in\HH(D_1)$. We recall that,
from the viewpoint of Stahl's theory, the first sheet of the two-sheeted R.s.\ associated
with the germ $f_\infty$ is the first sheet $\RS^{(1)}_2(w)\ni z^{(1)}$ of the R.s.\ $\RS_2(w)$
of the function $w^2=z^2-1$ defined by the inequality (see~\eqref{33})
$$
-g(z^{(1)})<-g(z^{(2)})
$$
with the canonical projection $\pi_2\colon\RS_2(w)\to\myo\CC$, $\pi_2(\zz)=z$,
where $\zz=(z,w)$ under the assumption that
$w/z\to1$ for $\zz\in\RS^{(1)}_2(w)$, $\zz\to\infty^{(1)}$, and
$w/z\to-1$ for $\zz\in\RS^{(2)}_2(w)$, $\zz\to\infty^{(2)}$.

From what has already been said above in \S\,\ref{s3s2} of this section it follows that,
in the case when at least one of the numbers~$A$ or~$B$ does not lie on the real line,
the Nuttall partition of the  R.s.\ $\RS_4(f)$ into sheets proceeds in a~nontrivial manner\,--\,namely,
not along closed curves on~$\RS_4(f)$ (corresponding, from the viewpoint of
the canonical projection $\pi_4\colon\RS_4(f)\to\myo\CC$, to the intervals
$\Delta_1$ and $\Delta_2$), but rather along the curves $\EE$, $\FF$ and~$\EE'$,
to which, for a~given canonical projection, there correspond some arcs $E$, $F$ and~$E'$
connecting pairwise the points $ -1,1$, $a,b$ and again the point $-1,1$. Correspondingly
the first sheet $\RS^{(1)}_4$ of the Nuttall partition into sheets of the R.s.\ $\RS_4(f)$
is now different from the first sheet of the Stahl surface
$\RS_2(w)$
(even though the active branch points $\zz=\pm1$ lying on the boundary of these two sheets are the same).
More precisely, we have
$\pi_2(\RS_2^{(1)}(w))=\myo\CC\setminus\Delta_1 \neq\myo\CC\setminus
E=\pi_4(\RS^{(1)}_4)$, because in a~general case  $E\neq\Delta_1$.
Additionally, as is easily seen from the results of \S\,\ref{s3s2}, we have $E\cap
F=\varnothing$, $F\cap E'=\varnothing$ ($E=E'$ in the real case).
These facts are of uttermost importance in the use of Hermite--Pad\'e polynomials in applications that
will be discussed in the concluding \S~\ref{s5} of the present paper (see Figs.~\ref{fig3.2}--\ref{fig3.5}).

Let $R_{\nn_1},R_{\nn_2}$ and $R_{\nn_3}$ be the remainder functions for the three-dimensional
multi\-indexes\footnote{In this section, the notation for
multiindices is different from that used before in~\S~\ref{s2}.}
$\nn_1:=(n,n-1,n-1)$, $\nn_2:=(n,n,n-1)$ and $\nn_3:=(n,n,n)$, respectively
(see \eqref{22},~\eqref{23} and~\eqref{18}). For $A,B\in\RR$, all these
multiindices are normal (for the definition of a~normal multiindex,
see~\cite{NiSo88}). This result is a~corollary of the fact of~\cite{FiLo11} that an arbitrary Nikishin system
is perfect (see also \cite{NiSo88} and~\cite{FiLoLoSo10}).
Therefore, we have as $z\to\infty$
\begin{equation}
R_{\nn_1}(z)=\frac{C_{\nn_1}}{z^{3n+1}}+\dots, \quad
R_{\nn_2}(z)=\frac{C_{\nn_2}}{z^{3n+2}}+\dots, \quad
R_{\nn_3}(z)=\frac{C_{\nn_1}}{z^{3n+3}}+\dots,
\label{45}
\end{equation}
where  $C_{\nn_1},C_{\nn_2},C_{\nn_3}\neq0$. For complex $A$ and $B$
this fact ceases to be true.

All three functions $R_{\nn_1},R_{\nn_2}$ and $R_{\nn_3}$ are (singlevalued)
meromorphic functions on the R.s.\ $\RS_4(f)$. One can easily find their divisors (of zeros and
poles). Namely, the following explicit representations hold:
\begin{equation}
\begin{aligned}
(R_{\nn_1})&=
(3n+1)\infty^{(1)}-n\infty^{(2)}-n\infty^{(3)}-n\infty^{(4)}-3(\aa+\bb)
+\sum_{j=1}^5\pp_{1,j}(n),\\
(R_{\nn_2})&=
(3n+2)\infty^{(1)}-n\infty^{(2)}-n\infty^{(3)}-n\infty^{(4)}-3(\aa+\bb)
+\sum_{j=1}^4\pp_{2,j}(n),\\
(R_{\nn_3})&=
(3n+3)\infty^{(1)}-n\infty^{(2)}-n\infty^{(3)}-n\infty^{(4)}-3(\aa+\bb)
+\sum_{j=1}^3\pp_{3,j}(n).
\end{aligned}
\label{46}
\end{equation}
In the real case, since a~Nikishin system is perfect, there can be no
cancelations (of zeros and poles) in~\eqref{46}; i.e., all
$\pp_{k,j}(n)\neq\aa,\bb$. Besides, it follows directly from~\eqref{45}
that all $\pp_{k,j}(n)\neq\infty^{(1)}$. Moreover, in can be shown
that, for any $n\in\NN$,
\begin{equation}
\pp_{k,j}(n)\in\RS_4(f)\setminus\myo{\RS^{(1)}_4}.
\label{46.2}
\end{equation}
However, this is not so in the complex case. It might happen that some of the above results also
hold in the case of a~``general position'', but this requires special consideration and
as yet remains unjustified. Unfortunately, by now the theory of Hermite--Pad\'e polynomials
is not capable of providing meaningful and pretty general results on the asymptotic behaviour of such
polynomials in the complex case.
Many of particular results obtained in the complex case depend, to a~greater or lesser degree,
on the possibility of reduction of this particular problem to the real setting
(see, for example,~\cite{BaGeLo18},~\cite{LoDi18}). Due to the imperfections of the presently developed
methods in the theory of Hermite--Pad\'e polynomials, practically all fairly general and well-known conjectures
on the asymptotic behaviour of such polynomials remain unsolved for a~long time (for various conjectures
in this direction, see  \cite{Nut84},~\cite{Sta88}, \cite{Apt08}, and also
\cite{Rak16}, \cite{Rak18},~\cite{Sue18},~\cite{Sue18b}). It is worth pointing out that
the conjecture on topological properties of
the canonical partition (presently called the Nuttall partition) of an arbitrary
Riemann surface into sheets posed by Nuttall \cite[Sect.~3]{Nut84} in 1984  was proved very recently in 2017 in~\cite{ChKoPaSu17}.

For the validity of the analogue of Theorem~\ref{th2} in the
complex case considered here, no additional constraints are required, because the
required result follows directly from a~more general recent result (see Theorem~2 of~\cite{ChKoPaSu17}).

\begin{theorem}\label{th5}
Let $Q_{n,j}$, $j=0,1,2,3$, be the Hermite--Pad\'e polynomials defined by~\eqref{18}.
Then, for any~$j$, the compact set
$E'=\pi_4(\EE')$ attracts as $n\to\infty$ all zeros of the polynomials
$Q_{n,j}$,
except  possibly a~finite  number  (independent of~$n$) of zeros.
Moreover,
\begin{equation}
\frac{Q_{n,2}}{Q_{n,3}}(z)\overset\mcap\longrightarrow -\bigl[f(z^{(1)})+f(z^{(2)})
+f(z^{(3)})\bigr],
\quad z\in\myo\CC\setminus E'.
\label{47}
\end{equation}
\end{theorem}

Let us now introduce the multiindices\footnote{In \S~\ref{s2} such multiindices were
denoted by $\nn_1',\nn_2'$ and $\nn_3'$, respectively. In the present section
many auxiliary notations are different from those introduced above.}
$\nn_1=(n,n-1,n-1)$, $\nn_2=(n-1,n,n-1)$ and $\nn_3=(n-1,n-1,n)$. Let $Q_{\nn_k,j}$,
$j=0,1,2,3$ be the corresponding type~I Hermite--Pad\'e polynomials for the tuple  $[1,f_\infty,f_\infty^2,f_\infty^3]$,
and $R_{\nn_k}(z)$, $k=1,2,3$, be the corresponding remainder functions
\begin{equation}
R_{\nn_k}(z)=(Q_{\nn_k,0}+Q_{\nn_k,1}f_\infty+Q_{\nn_k,2}f_\infty^2+Q_{\nn_k,3}
f^3_\infty)(z)
=O\(\frac1{z^{3n+1}}\),\quad z\to\infty.
\label{47.2}
\end{equation}
The functions $R_{\nn_k}$ are singlevalued meromorphic functions on the R.s.\ $\RS_4(f)$, which are
completely defined by its divisor,
\begin{equation}
(R_{\nn_k})=(3n+1)\infty^{(1)}-n\infty^{(2)}-n\infty^{(3)}
-n\infty^{(4)}-3(\aa+\bb)+\sum_{j=1}^5\pp_{k,j}(n).
\label{48}
\end{equation}

Given two arbitrary distinct points $\zz_1$ and $\zz_2$,
$\zz_1,\zz_2\in\RS_4(f)$, we denote by $G(\zz_1,\zz_2;\zz)$
an Abelian integral of the third kind on the R.s.\ $\RS_4(f)$ with purely imaginary periods and
logarithmic singularities only at the points $\zz_1$ and $\zz_2$, which in the
corresponding  local coordinate read as
\begin{equation}
\begin{aligned}
G(\zz_1,\zz_2;\zz)&=-\log{\zeta}+O(1),\quad \zz\to\zz_1,\\
G(\zz_1,\zz_2;\zz)&=\log{\zeta}+O(1),\quad \zz\to\zz_2.
\end{aligned}
\label{49}
\end{equation}
Consequently,
\begin{equation}
\Psi(\zz_1,\zz_2;\zz):=\exp\bigl\{G(\zz_1,\zz_2;\zz)\bigr\}
\label{50}
\end{equation}
is a multivalued analytic function on the R.s.\ $\RS_4(f)$ with singlevalued modulus,
a~pole at the point $\zz=\zz_1$, and a~zero at the point $\zz=\zz_2$. According to \eqref{48}
and~\eqref{50}, for any $k=1,2,3$ the remainder function $R_{\nn_k}$ can be written as
\begin{equation}
R_{\nn_k}=C_k\exp\biggl\{n\sum\limits_{j=2}^4G(\infty^{(j)}),\infty^{(1)};\zz)
+\sum\limits_{j=1}^5G(\aa_j,\pp_{k,j}(n);\zz)+G(\aa_6,\infty^{(1)};\zz)
\biggr\},
\label{51}
\end{equation}
where $C_1,C_2,C_3\neq0$ are the normalized constants and for convenience we have denoted
the poles of the function $f^3$ on the R.s.\ $\RS_4(f)$ by
$\{\aa_1,\dots,\aa_6\}$, where each pole appears as many times as its multiplicity.

To prove analogues of Theorems~\ref{th2} and~\ref{th3} we shall require some additional conditions.
As was already pointed out, these conditions, which are not related to the essence of the matter,
appear only due to the imperfections of the presently developed methods for
the study of asymptotic properties of Hermite--Pad\'e polynomials.
It is also worth noting that the imposed additional conditions are similar to \eqref{46.2},
but it seems that they cannot be expressed directly in terms of the original germ $f_\infty$.

\begin{theorem}\label{th6}
Assume that for some infinite sequence $\Lambda\subset\NN$
there exists a~neigh\-bour\-hood~$\bU$ of the point  $\zz=\infty^{(1)}$ such that all points
$\pp_{k,j}(n)$, as defined by~\eqref{46}, lie outside this neighbourhood,
$\pp_{k,j}(n)\in\RS_4(f)\setminus\myo{\bU}$. Then the following assertions hold:

1) the compact set $F=\pi_4(\FF)$
attracts as $n\to\infty$, $n\in\Lambda$, all zeros of the polynomials
$S_{2n,1}:=Q_{\nn_1,1}Q_{\nn_2,3}-Q_{\nn_1,3}Q_{\nn_2,1}$ and
$S_{2n,2}:=Q_{\nn_1,2}Q_{\nn_2,3}-Q_{\nn_2,2}Q_{\nn_1,3}$ except  possibly a~finite  number  (independent of~$n$) of zeros;

2) if  $n\to\infty$, $n\in\Lambda$, then
\begin{equation}
\frac{S_{2n,1}}{S_{2n,2}}(z)\overset\mcap\longrightarrow-\bigl[f(z^{(1)})+f(z^{(2)})\bigr],
\quad z\in\myo\CC\setminus F.
\label{48.2}
\end{equation}
\end{theorem}

The proof of Theorem~\ref{th6} is based, first, on the determinant representation of the
polynomials  $S_{2n,j}$ (see \eqref{24}), which in the present paper we used as a~definition
of these polynomials, and second, on an explicit representation of the remainder function~\eqref{51}
in terms of Abelian integrals of the third kind with purely imaginary periods. These
representations follows from \eqref{46} and from explicit formulae for
type~I Hermite--Pad\'e polynomials
which follow from such a~representation of the remainder functions
(see \cite{Nut84}, \cite[formula~(38)]{ChKoPaSu17}).

Since  due to the lack of general methods capable of dealing with asymptotic properties
of type~II Hermite--Pad\'e polynomials in the complex case which are based directly on
definition~\eqref{11}--\eqref{13} of these polynomials, we use
here the determinant representation~\eqref{28}; the singularity of some determinants
related to the explicit representation of the remainder functions~\eqref{51} will be of key importance for us.
Namely, we set
\begin{equation}
\Psi_{n,k}(\zz):=\exp\biggl\{
\sum_{j=1}^5G(\aa_j,\pp_{k,j}(n);\zz)\biggr\},\quad k=1,2,3.
\label{52}
\end{equation}
The functions $\Psi_{n,k}$ are multivalued analytic functions on the R.s.\ $\RS_4(f)$  with
singlevalued modulus and multivaluedness of the  same character\,--\,all three functions become
singlevalued after multiplication by the same factor
$$
\exp\biggl\{\sum_{k=2}^4G(\infty^{(k)},\infty^{(1)};\zz)+G(\aa_6,\infty^{(1)};\zz)
\biggr\}.
$$
We set
\begin{equation}
\Delta_n(\zz):=\mdet\bigl[\Psi_{n,k}(z^{(j)})\bigr]_{j=1,2,3;k=1,2,3}.
\label{53}
\end{equation}

The following result holds.

\begin{theorem}\label{th7}
Assume that there exists an infinite sequence $\Lambda\subset\NN$ such that,
for some finite set of points
$\bM:=\{\zz_1,\dots,\zz_M\}\subset\RS_4(f)$ independent of $n\in\Lambda$ and
an arbitrary compact set
$\KK\subset\RS_4(f)\setminus\bM$,
\begin{equation}
0<C_1(\KK)\leq|\Delta_n(\zz)|\leq C_2(\KK)<\infty,\quad \zz\in\KK,
\quad n\in\Lambda.
\label{54}
\end{equation}
Then the following assertions hold:

1) if $n\to\infty$, $n\in\Lambda$, then the compact set $E=\pi_4(\EE)$ attracts
all zeros of the polynomials $Q_{3n}(z)$, except  possibly a~finite  number  (independent of~$n$) of zeros.

2) if $n\to\infty$, $n\in\Lambda$, then
\begin{equation}
\frac{P_{3n,1}}{Q_{3n}}(z)\overset\mcap\longrightarrow f(z^{(1)}),\quad z\in\myo\CC\setminus
E.
\label{56}
\end{equation}
\end{theorem}

As in the real case, relations \eqref{47},~\eqref{48.2}
and~\eqref{56} guarantee the possibility to recover the function $f(z)$ in the domain
$\RS_4(f)\setminus\myo{\RS^{(4)}_4}$ by successive recoveries of the values
$f(z^{(1)}),f(z^{(2)})$ and $f(z^{(3)})$ of this function
on the first, second, and third sheets of the R.s.\ $\RS_4(f)$
ordered in accordance to the Nuttall partition.

\begin{remark}\label{rem4}
It is worth pointing out that the Nuttall partition is related to
a~selected point at which a~germ of the original function~$f$ is
defined. Here we consider a~germ $f_\infty$ defined at the point
$z=\infty$. Both the Pad\'e polynomials and all Hermite--Pad\'e
polynomials considered here were defined in accordance with this germ
at the point $z=\infty$. However, if the original germ is given at some
other point $z_0\in\CC$, then all the corresponding constructions of
Pad\'e polynomials and Hermite--Pad\'e polynomials should be adopted
appropriately to the point $z_0$. Now the Nuttall partition is given
with the help of an Abelian integral with purely imaginary periods and
having suitable logarithmic singularities of the form \eqref{42.2} at
points of the set $\pi_4^{-1}(z_0)$, as written in the local
coordinate~$\zeta$.
\end{remark}

\section{Concluding remarks}\label{s4}

\subsection{}\label{s4s1}

Let us summarize to some extent what has been said above.

Given a function $f\in\SZ$ of the form \eqref{44}, we have the following results (here
we retain the above notation, but for convenience we introduce new labeling). If $n\to\infty$, then
\begin{align}
\frac{P_{3n,1}}{Q_{3n}}(z)&\overset\mcap\longrightarrow f(z^{(1)}),\quad
z\in\myo\CC\setminus E;
\label{91}\\
\frac{S_{2n,1}}{S_{2n,2}}(z)&\overset\mcap\longrightarrow-\bigl[f(z^{(1)})+f(z^{(2)})\bigr],
\quad z\in\myo\CC\setminus F;
\label{92}\\
\frac{Q_{n,2}}{Q_{n,3}}(z)&\overset\mcap\longrightarrow -\bigl[f(z^{(1)})+f(z^{(2)})
+f(z^{(3)})\bigr],
\quad z\in\myo\CC\setminus E'.
\label{93}
\end{align}
For real parameters $A$ and $B$, these results hold as $n\to\infty$
for all $n\in\NN$; for complex parameters $A$ and $B$ these results are known at present to be
true only along some subsequences.
Clearly, using the limit relations~\eqref{91}--\eqref{93} it is possible,
when considering them as a~whole, to recover in succession the values of the
function $f(z)$ at the three Nuttall sheets of its R.s.\ $\RS_4(f)$, and hence,
in the domain $\RS_4(f)\setminus\myo{\RS^{(4)}_4}$ lying on the R.s.\ $\RS_4(f)$.
As has already been mentioned, on the last fourth sheet it is not possible to recover
the function~$f$ by the above procedure.
Nevertheless, it follows from \eqref{91}--\eqref{93} that our new approach to the
solution of the problem of efficient analytic continuation of a~given germ is vastly superior to
the approach based on Pad\'e approximants. Indeed, using this approach one can
extend the original germ to the {\it domain~$\RS_4(f)\setminus\myo{\RS^{(4)}_4}$
lying on the Riemann surface}
of the multivalued analytic function corresponding to this germ. At the same time,
with the help of Pad\'e approximants, the analytic continuation can be carried out only to the  {\it domain
lying on the Riemann sphere}. Note that in some sense the domain $\RS_4(f)\setminus\myo{\RS^{(4)}_4}$
is a~three-sheeted covering of the Stahl domain (see~Fig.~\ref{fig}).

The theoretical results presented here have a~very special character, but
in spite of this they are fairly typical to demonstrate the advantages of the
new approach.

We note once again that for any $n\in\NN$ the number of Laurent coefficients
involved in the construction consisting of
three rational functions under consideration  (see \eqref{91}--\eqref{93}) is at most $3n+4$.

\subsection{}\label{s4s2}
It would be natural to compare the approach proposed here not only with the diagonal PA, but also with
other available methods of analytic continuation.

\subsubsection{Universal matrix method of summation of power series}\label{s4s2s1}
An account of these methods may be found in \cite{Ara84},~\cite{ArLu02},~\cite{Ara03} and
in the references given therein.

The most well known of these methods is the Mittag-Leffler method, which
we now  briefly discuss. This a~linear method of extension
(summation) of an element~$f_{z_0}$ to the Mittag-Leffler star.

A~domain $G\subset\myo\CC$ is called {\it star-shaped} with respect to a~point
$z_0\in G$ if, for any point $\zeta\in G$, the interval
$[z_0,\zeta]$ also lies in~$G$.

\begin{theoremML}
There exists a~table of numbers
$\{b^{(n)}_0,b^{(n)}_1,\dots,\allowbreak b^{(n)}_{k_n}\}$, $n=0,1,\dots$,
such that the following property holds for any domain~$G$ star-shaped with respect to some point $z_0\in G$
and for any function~$f$ holomorphic in the domain~$G$.

Assume that a function $f\in\HH(G)$ is given by a~power series\footnote{%
Only in this section \ref{s4s2s1} we follow the tradition of the book \cite{Bib67}
to use the same letter to denote the Taylor coefficients of an analytic function
and the function itself.}
at the point $z=z_0$
\begin{equation}
\label{e1}
f(z)=\sum_{\nu=0}^\infty f_\nu(z-z_0)^\nu,\qquad |z-z_0|<R_0.
\end{equation}
Then
\begin{equation}
\label{e2}
f(z)=\sum_{n=0}^\infty \Bigl\{b^{(n)}_0f_0+b^{(n)}_1f_1(z-z_0)
+\dotsb+b^{(n)}_{k_n}f_{k_n}(z-z_0)^{k_n}\Bigr\},
\quad z\in G,
\end{equation}
uniformly inside the domain~$G$.
\end{theoremML}

\begin{remark}\label{rem5}
Assume that at a point $z=z_0$ we are given two power series
$g=\sum\limits_{\nu=0}^\infty g_\nu(z-z_0)^\nu$ and
$d=\sum\limits_{\nu=0}^\infty d_\nu(z-z_0)^\nu$. Then the series
\begin{equation}
\label{ad1}
h(z):=(g\circ d)(z)=\sum_{\nu=0}^\infty g_\nu d_\nu (z-z_0)^\nu
\end{equation}
is called the {\it Hadamard composition} of the series $g$ and~$d$; see~\cite{Bib67}.
The expression in the curly brackets under the summation sign in~\eqref{e2} is a~Hadamard composition
$f\circ P_n$ of the original function $f(z)$ and the polynomial
$P_n(z-z_0):=b^{(n)}_0+b^{(n)}_1(z-z_0)+\dots b^{(n)}_{k_n}(z-z_0)^{k_n}$. So, \eqref{e2}~can
be written as
\begin{equation}
\label{ad2}
f(z)=\sum_{n=0}^\infty (f\circ P_n)(z-z_0),
\quad z\in G.
\end{equation}
\end{remark}

We have thus  obtained (see \eqref{e2}) an expansion of the function $f(z)$ in a~series of polynomials
uniformly convergent on the set~$K$, and therefore, inside~$G$, since  $K$~is arbitrary.
The coefficients $b^{(n)}_0,b^{(n)}_1,\dots,b^{(n)}_{k_n}$, $n=0,1,2,\dots$,
are independent of~$f$ and can be evaluated once for all. In addition to them,
the formula involves only the coefficients $\{f_\nu,\nu=0,1,2,\dotso\}$, i.e., the
coefficients of the power series~\eqref{e1} defining the analytic element
$f_{z_0}=\bigl\{f(z),|z-z_0|<R_0\}$ that extends to the entire domain~$G$ as
a~holomorphic (singlevalued analytic) function; see \eqref{e2}. This expansion
(known as the Mittag-Leffler expansion) clearly solves the extension problem of
the element $f_{z_0}$ {\it along straight rays} (see \cite[Ch.~8,
\S~5.5.3]{Mar68}).

In~\cite[Ch.~8, \S~5.5.4]{Mar68} a~suitable sequence of polynomials
$\{P_n(w)\}$ was constructed explicitly. This construction is based on
the method proposed by Painlev\'e \cite{Pai05}.

The Mittag-Leffler method is a~{\it linear method}
of summation of power series in the Mittag-Leffler star.
The Mittag-Leffler star itself cannot be recovered using the Mittag-Leffler representation\footnote{%
At least the authors are unaware  on any results in this direction.}.
Thus, as in the Stahl theory, summation of power series takes place in some {\it
domain lying on the Riemann sphere~$\myo\CC$.}

\subsubsection{Shafer quadratic approximants}\label{s4s2s2}
For a tuple $[1,f_\infty,f_\infty^2]$ and an arbitrary $n\in\NN$, the {\it Shafer approximants} (see \cite{Sha74})
are defined in terms of type~I Hermite--Pad\'e polynomials for
the two-dimensional multiindex $\nn:=(n,n)$ (see~\eqref{81}) as follows:
\begin{equation}
\Sh_{n;1,2}(z):=\frac{-Q_{\nn,1}(z)\pm\bigl[(Q_{\nn,1}^2-4Q_{\nn,0}Q_{\nn,2})(z)
\bigr]^{1/2}}{2Q_{\nn,2}(z)}.
\label{94}
\end{equation}
Sometimes Shafer approximants are called Hermite--Pad\'e approximants.
The difficulties with Shafer quadratic approximants
in finding the analytic continuation of power series are obvious from
definition~\eqref{94}. First, this is the necessity of taking the root of the discriminant
$D_n(z):=(Q_{\nn,1}^2-4Q_{\nn,0}Q_{\nn,2})(z)$, which is a~polynomial of degree~$2n$
having $2n$ (complex) roots. Second, the sign
``$+$'' or ``$-$'' of the square root in~\eqref{94} should be chosen appropriately.
After all, Shafer approximants are clearly capable of delivering in principle
the values of a~multivalued analytic function only on the first two sheets of the
three-sheeted Nuttall surface associated with this function. Clearly, in view of
\eqref{83} and~\eqref{84} (see \S\,\ref{s4s3}
below), the type~II and~I Hermite--Pad\'e polynomials with the two-dimensional multiindex $\nn:=(n,n)$
are much better in coping with this problem.

Some, seemingly first, theoretical results on the limit distribution of the zeros of the discriminant $D_n$ and
on the convergence of Shafer approximants can be found in~\cite{KoKrPaSu16}.

\subsection{}\label{s4s3}
The Weierstrass approach leads to the following natural problem: find the value of the Weierstrass extension
of the element $\myt{f}_0$ along some curve without reexpansion of power series.
From the results of the present paper it follows that this can surely be achieved with the
help of Hermite--Pad\'e polynomials associated with multiindices of large dimension.

The results formulated in the present paper (see Theorems~\ref{th1}--\ref{th7}) will be
used for the numerical solution of the above problem in Example~\ref{ex3}.

Without going into further details we show here how the problem associated with the construction of
the analytic continuation using PA can be circumvented by invoking
type~I and~II polynomials defined similarly to \eqref{11}--\eqref{13} and~\eqref{18},
but, respectively, for the tuple
$[1,f_\infty,f_\infty^2]$ and the pair of germs $f_\infty,f_\infty^2$ (note that in this case
the type~I and~II Hermite--Pad\'e polynomials can be looked upon as two variants of a~natural extension of
Pad\'e polynomials). Such polynomials are well known (see, first of all, \cite{Nut84}, and also~\cite{NiSo88}).

We recall the following fact, which is a fairly special case of one general fact resulting from Stahl's theory.

For an arbitrary function $f\in\SZ$, the Stahl compact set is the closed interval $[-1,1]$.
The Riemann surface of the function $f\in\SZ$ has at least four sheets.
But in accordance with Stahl's theory, with an arbitrary function $f\in\SZ$ one can
{\it associate} in a~certain sense the two-sheeted R.s.\ $\RS_2(w)$ of the function~$w$,
where $w^2=z^2-1$. Such a~R.s.\ will be called the {\it Stahl surface} of a~function $f\in\SZ$ or
the {\it Stahl surface associated with the  function $f\in\SZ$}.
The fact this definition is quite meaningful follows from the paper \cite{ApYa15} by Aptekarev and Yattselev,
which proves that the problem of strong asymptotics of Pad\'e polynomials
can be solved in a~fairly wide class of multivalued analytic functions
exactly in the terms associated with the two-sheeted Stahl surface (see also~\cite{Nut84}). Note that in general
the two-sheeted Stahl surface is a~hyperelliptic R.s.

The original  germ  $f_\infty\in\HH(\infty)$ extends to the first sheet
of the Stahl surface (i.e., the R.s.\ $\RS_2(w)$ of the function~$w$) as a~singlevalued
holomorphic function. But for a~singlevalued extension of this germ
to the second sheet of this surface
suitable cuts should be drawn on this sheet. On the entire first
(open) sheet of the Stahl surface the function $f(z)$ is recovered using
diagonal PA.

In some cases, for a function $f\in\SZ$ (see \cite{IkKoSu18}
and cf.~\cite{RaSu13}) it proves possible to show that with a~germ $f_\infty\in\HH(\infty)$
of~$f$ one can associate a~{\it three-sheeted} Riemann surface $\RS_3=\RS_3(w)$
defined by some irreducible cubic equation (see \eqref{85}) and such that, for
the Nuttall partition of this surface into sheets
(see \cite{Nut84},~\cite{RaSu13}, \cite{MaRaSu16}, \cite{ChKoPaSu17}) and
when the points $z=\infty\in\myo\CC$ and
$\zz=\infty^{(1)}\in\RS_3(w)$ are identified, the original germ $f_\infty$ of the function $f\in\SZ$
extends to the domain $\RS_3(w)\setminus\myo{\RS^{(3)}_3}$ as a~singlevalued
meromorphic  function. For a~singlevalued extension of this germ to the third sheet $\RS^{(3)}_3$ of this R.s.\ $\RS_3(w)$
one should draw appropriate cuts on this sheet.
In contrast to Stahl's theory, so far the existence of such a~three-sheeted R.s.\
was established only for rather special cases of multivalued analytic functions.
Following~\cite{RaSu13} and~\cite{ChKoPaSu17}, such a~three-sheeted R.s.\ will be called
a~{\it Nuttall surface associated with the germ $f_\infty$}, or simply,
a~{\it Nuttall surface}.

Let $\pi_3\colon\RS_3(w)\to\myo\CC$ be the canonical projection,
$F^{(1,2)}$ be the boundary between the first and second sheets, $F^{(2,3)}$ be the boundary
between the second and third sheets of the Nuttall partition of the R.s.\ $\RS_3(w)$.
We set $E_2:=\pi_3(F^{(1,2)})$, $F_2:=\pi_3(F^{(2,3)})$.

For a germ $f_\infty\in\HH(\infty)$, the tuple of germs\footnote{Here and in what follows
we assume that the three functions $1,f,f^2$ are independent over the field $\CC(z)$.}
$[1,f_\infty,f_\infty^2]$, an arbitrary $n\in\NN$, and the two-dimensional multiindex
$\nn:=(n,n)$, we define (not uniquely) the type~I Hermite--Pad\'e polynomials
$Q_{\nn,0}, Q_{\nn,1}$ and $Q_{\nn,2}$, $Q_{\nn,j}\in\PP_n$, $Q_{\nn,j}\not\equiv0$, by the relation
\begin{equation}
(Q_{\nn,0}+Q_{\nn,1}f_\infty+Q_{\nn,2}f_\infty^2)(z)=O\(\frac1{z^{2n+2}}\),
\quad z\to\infty.
\label{81}
\end{equation}
For a~pair of germs $f_\infty,f_\infty^2$ and $n\in\NN$, the type~II Hermite--Pad\'e polynomials
$Q_{2n},P_{2n,1}$ and $P_{2n,2}$, $Q_{2n},P_{2n,j}\in\PP_{2n}$,
$Q_{2n},P_{2n,j}\not\equiv0$,
are defined (not uniquely) by the relations
\begin{equation}
\begin{aligned}
(Q_{2n}f_\infty-P_{2n,1})(z)&=O\(\frac1{z^{n+1}}\),\quad z\to\infty,\\
(Q_{2n}f_\infty^2-P_{2n,2})(z)&=O\(\frac1{z^{n+1}}\),\quad z\to\infty.
\end{aligned}
\label{82}
\end{equation}
In some cases it proves possible  (see \cite{IkKoSu18},~\cite{LoKoSu18})
to show that as $n\to\infty$
\begin{align}
\frac{P_{2n,1}}{Q_{2n}}(z)&\overset\mcap\longrightarrow f(z^{(1)}),
\quad z\in\myo\CC\setminus{E_2},
\label{83}\\
\frac{Q_{\nn,1}}{Q_{\nn,2}}(z)
&\overset\mcap\longrightarrow -\bigl[f(z^{(1)})+f(z^{(2)})\bigr],
\quad z\in\myo\CC\setminus{F_2}.
\label{84}
\end{align}
Note that if $n\to\infty$, then the compact set $E_2$ attracts ``almost all''\footnote{As before, here we speak about the
limit distribution of zeros, which means that $o(n)$ many zeros
remain uncontrolled as $n\to\infty$.}
zeros of the polynomials $Q_{2n}$ and $P_{2n,1}$
and  the compact set $F_2$ attracts ``almost all'' zeros of the polynomials $Q_{\nn,1}$ and
$Q_{\nn,2}$. Similarly to \eqref{17} and~\eqref{21}, the rate of convergence in~\eqref{83}
and~\eqref{84} can be characterized in terms related to the scalar potential theory equilibrium problem.
From relations~\eqref{83} and~\eqref{84} one can recover in succession
the values of the function $f_\infty\in\MM(\RS_3(w)\setminus\myo{\RS_3^{(3)}})$ on the first two sheets
of the Nuttall surface. It is worth pointing out that in \eqref{83} and~\eqref{84}
we mean the Nuttall partition into sheets of the R.s.\ $\RS_3=\RS_3(w)$, where the function~$w$
satisfies the irreducible third-degree equation
\begin{equation}
w^3+r_2(z)w^2+r_1(z)w+r_0(z)=0,
\label{85}
\end{equation}
here $r_0,r_1,r_2\in\CC(z)$ are some rational functions of~$z$ (cf.\ the equation
$w^2-(z^2-1)=0$, which defines the Stahl surface for an arbitrary function
$f\in\SZ$).
In a~manner completely similar to \eqref{17} and~\eqref{21}, the
rate of convergence in \eqref{83} and~\eqref{84} can be characterized
in terms of some scalar potential theory equilibrium problem
(see \cite{RaSu13},~\cite{MaRaSu16},~\cite{LoKoSu18} and
cf.\ \cite{ApBoYa17},~\cite{LoVa18}).

\begin{remark}\label{rem6}
Note that $E_2\cap F_2=\varnothing$ and, generally speaking, the compact set $E_2$
is different from the above compact set~$E$, and  $F_2$~is different from the analogous
compact set~$F$ (see \eqref{48.2} and~\eqref{56}). More precisely,
$\pi_3(F^{(1,2)})=E_2\neq E=\pi_4(\Gamma^{(1,2)})$ and
$\pi_3(F^{(2,3)})=F_2\neq F=\pi_4(\Gamma^{(2,3)})$ (the coincidence may take place only in exceptional cases,
for example, for real  $A$ and~$B$ in representation~\eqref{44}). This fact has great value
in the applications of our theoretical results that we shall discuss below;
see, for instance, Example~\ref{ex3} and Figs.~\ref{fig3.2}--\ref{fig3.5}.
\end{remark}

\begin{remark}\label{rem7}
In definitions~\eqref{81} and~\eqref{82} the following changes are required in the case when the
original germ of an analytic function~$f$ is given not
at the point $z=\infty$, but at some point $z_0\in\CC$, $f_{z_0}\in\HH(z_0)$.
The Nuttall partition into
sheets of the three-sheeted R.s.\ associated with this germ is now defined with respect to the
highlighted point~$z_0$.

In what follows, when discussing possible applications of the method proposed here, we shall
consider the germs given at the point $z=0$.
\end{remark}

\section{Some possible applications}\label{s5}

\subsection{}\label{s5s1}
In this section we shall discuss some possible of applications of the above theoretical results
in numerical mathematics for approximate solution of problems
based on the analytic continuation of the original data
(for example, the analytic continuation with respect to the small parameter).
We give three examples based on functions of the form \eqref{44}, in which
the quantities $A$ and~$B$ have the same imaginary part and opposite real parts;
i.e., \begin{equation}
f(z):=\[\(A_1+iA_2-\frac1{\pfi(z)}\)\(-A_1+iA_2-\frac1{\pfi(z)}\)\]^{-1/2},
\label{61}
\end{equation}
where $A_1,A_2\in\RR\setminus\{0\}$.
Here, from the variable~$z$ we change to the variable $\myt{z}$, where
$z=(\myt{z}-a)i/b$ and $a,b\in\RR\setminus\{0\}$ (clearly,
such a~transform rotates the interval $[-1,1]$ clockwise by the angle  $\pi/2$ with
magnification by~$|b|$ and translation by $|a|$ of the rotated and expanded interval to the right or left half-plane
depending on the sign of~$a$). Correspondingly,
instead of~$\pm1$ in the plane $\CC_{\myt{z}}$ there appears a~pair of complex conjugate
second-order branch points  $\myt{z}_1=a\pm ib$. The second pair of branch points (also of order two)
corresponds to the choice of the other branch
of the inverse of the Zhukovskii function; it also consists of a~pair of complex conjugate points
$\myt{z}_{1,2}=\myt{z}_{1,2}(a,b)$ lying outside the interval $[a-ib,a+ib]$ and
corresponding to $A_1+iA_2$ and $-A_1+iA_2$.
In Examples~\ref{ex1}--\ref{ex3} we choose concrete  values of
$A_1,A_2,a$ and~$b$. In the first case, $a<0$, and in the second and third cases,
$a>0$ ($|a|$ remains unchanged).
Concrete values of $A_1,A_2,a$ and $b$ are immaterial, and hence are not given.
At last, we change the variable $\zeta=1/\myt{z}$, transforming the point
$\myt{z}=\infty$ into the point $\zeta=0$, the
Stahl compact set and the compact sets $E_2$, $F_2$, $E,F$ and $E'$ are transformed accordingly (see~\S~\ref{s3},
Remark~\ref{rem4}, and \S~\ref{s4}, Remark~\ref{rem7}). In particular,
the Stahl compact set, which is the interval $[a-ib,a+ib]$
in the $\myt{z}$-plane, is moved into a~circular arc passing
through the point $\zeta=0$ and the points $\zeta=1/(a-ib)$, $\zeta=1/(a+ib)$.
Thus, in view of~\eqref{61}, the function considered in the present section reads as
\begin{align}
f(z)=\myt{f}(\zeta)=
&
\[\(A_1+iA_2+\frac{i\zeta}{(1-a\zeta)/b+\sqrt{((1-a\zeta)/b)^2
+\zeta^2}}\)
\right.
\notag\\
&\times
\left.\(-A_1+iA_2+\frac{i\zeta}{(1-a\zeta)/b+\sqrt{((1-a\zeta)/b)^2
+\zeta^2}}\)
\]^{-1/2}
\label{62}
\end{align}
Let us explain which properties of Hermite--Pad\'e polynomials will be demonstrated on these three
examples.

\begin{figure}[!ht]
\centerline{
\includegraphics[width=15cm,height=15cm]{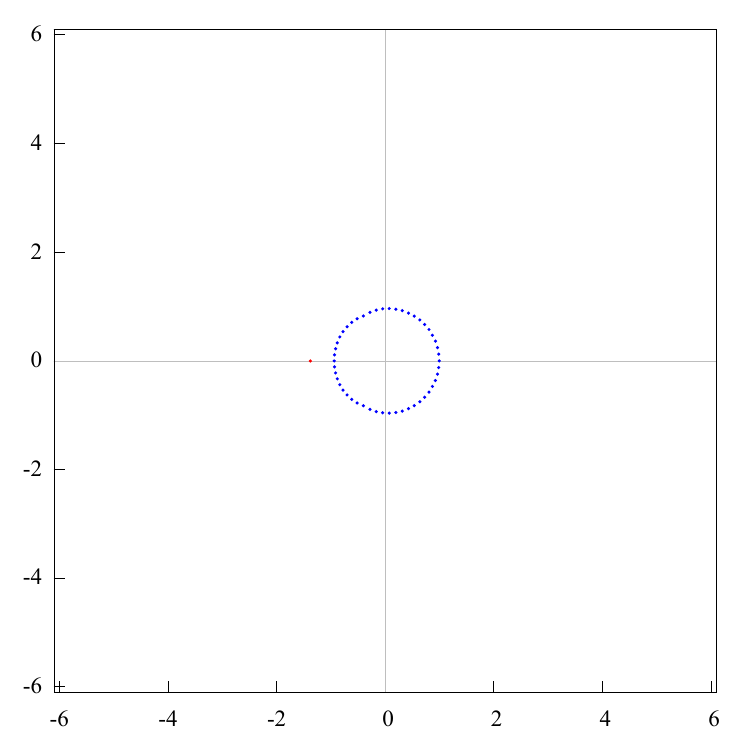}}
\vskip-6mm
\caption{Example~\ref{ex1}.
The zeros (dark blue points) of the partial sum $S_{50}$
for the function $\myt{f}$ from Example~\ref{ex1}.
In complete accord with the classical Jentzsch--Szeg\H o's theorem
(see \cite{Jen16},~\cite{Sze21}), these zeros model
the circle (the boundary of the convergence disc of the power series $\myt{f}_0\in\HH(0)$).
}
\label{fig00}
\end{figure}

\begin{figure}[!ht]
\centerline{
\includegraphics[width=15cm,height=15cm]{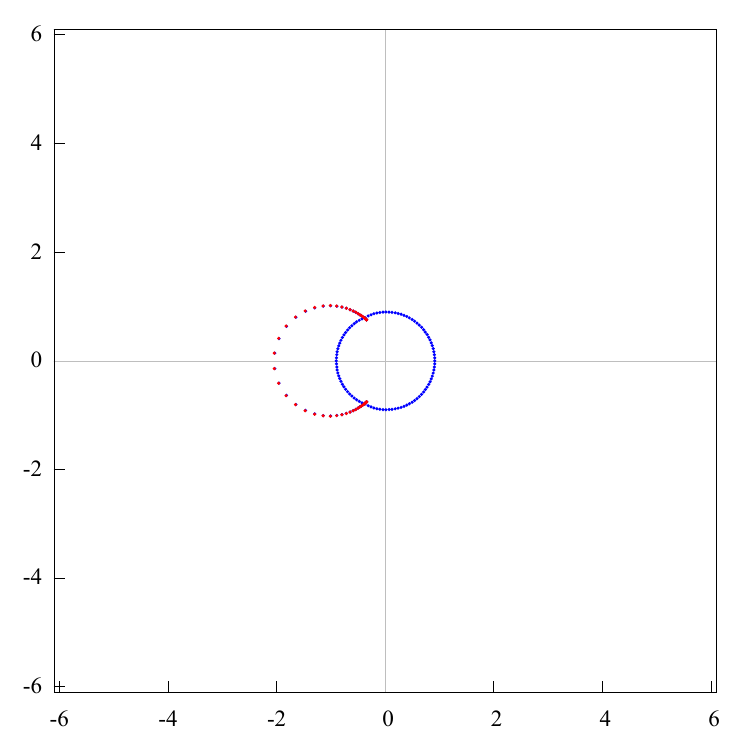}}
\vskip-6mm
\caption{Example~\ref{ex1}.
The zeros of the partial sum $S_{100}$ (dark blue points) and the poles of the diagonal
PA $[50/50]_{\myt{f}_0}$ (red points).
It is clearly seen that from the zeros and poles
of PA one can find the branch points of a~multivalued analytic function,
while this cannot be achieved using the zeros of the partial sums. The convergence radius
of power series can be evaluated more accurately using the diagonal PA, rather than by employing
the partial sums. The same number of
Taylor coefficients of this series (101 coefficients $c_0,c_1,\dots,c_{100}$)
were used to find the partial sum $S_{100}$ and the diagonal PA $[50/50]_{\myt{f}_0}$ .
}
\label{fig01}
\end{figure}

So, an original function of the form \eqref{62} is defined by its series $\myt{f}_0$,
which converges near the point $\zeta=0$, $\myt{f}_0\in\HH(0)$. We are concerned with the
problem of Weierstrass analytic continuation of the given power series $\myt{f}_0$
from the point $x=0$ along the real line to the point $x=6$ (it will be clear later that
the choice of the point $x=6$ as an end-point of the extension is immaterial).
Since the point $x=6$ is outside of the disc of convergence of the power series
$\myt{f}_0$, the value
$\myt{f}(6)$ cannot be (approximately) evaluated using partial sums of
the series $\myt{f}_0$ of arbitrarily high order.
Of course, this fact is a~direct consequence of the Cauchy--Hadamard formula.
However, there is another more transparent method to illustrate this fact
in the framework of the approach discussed here. Namely, we analyze
the limit distribution of free zeros and poles of rational approximants.
Figure~\ref{fig00} depicts the numerically obtained distribution of the zeros of the partial sum
$S_{50}$ of the series $\myt{f}_0$. It is transparent that the zeros of the partial sum $S_{50}$
model the circle (the boundary of the convergence disc of the power series $\myt{f}_0$). This is in full accord
with the classical Jentzsch--Szeg\H o theorem
(see \cite{Jen16},~\cite{Sze21}), according to which, along some subsequence depending
on the function, the limit distribution of the zeros of partial sums exists and coincides with the Lebesgue measure of this circle.
Looking ahead, we point out that Fig.~\ref{fig01} depicts the zeros of the
partial sum $S_{100}$ and the zeros and poles of the diagonal PA $[50/50]_{\myt{f}_0}$ for the function from
Example~\ref{ex1}. It is transparent that using zeros and poles
of the diagonal PA one can find the branch points of a~multivalued analytic function,
but this cannot be done using the zeros of partial sums. Even the radius of convergence
of power series can be evaluated more accurately using the diagonal PA, rather than with the help of
partial sums. It should be pointed out here that the same number
of the coefficients of Taylor series
was used to calculate the partial sum $S_{100}$ and the diagonal PA $[50/50]_{\myt{f}_0}$
(namely, 101 coefficients $c_0,c_1,\dots,c_{100}$ were used).
Figure~\ref{fig00}, which is formally unrelated to PA and Hermite--Pad\'e polynomials,
completely fits the general idea underlying the use of constructive rational approximants:
the limit distribution  of their {\it free zeros and poles} is most immediately related to the boundary
of that domain in which these rational approximants
realize an efficient analytic continuation of the original  power series.
Indeed, all poles of the partial sums of the power series are fixed at the
infinity point. Hence the boundary of the convergence domain of partial sums
is controlled by the limit distribution of {\it free} zeros of these partial sums.

To be more precise, let us clarify what is meant by a~(singlevalued) analytic continuation
of the power series $\myt{f}_0$ from the point $x=0$ to the right along the real line
to the point $x=6$. Here we
speak about the Weierstrass approach to the concept of an analytic function, which leads to the concept of
a~{\it global analytic function}. This approach resides on the concept of an
{\it analytic element} (or simply an {\it element}) of an
analytic function, which is defined as the convergent power series defined at a~fixed point $\zeta_0$ and
having the radius of convergence $R_0>0$. This will be briefly written
as $\myt{f}_0\cong(\zeta_0,\DD(\zeta_0,R_0))$, where
$\DD(\zeta_0,R_0):=\{\zeta\in\CC:|\zeta-\zeta_0|<R_0\}$.
An element $\myt{f}_1\cong(\zeta_1,\DD(\zeta_1,R_1))$ is called the  {\it
direct} analytic continuation of an $\myt{f}_0$ if $D:=\DD(\zeta_0,R_0)\cap
\DD(\zeta_1,R_1)\neq\varnothing$ and $\myt{f}_0(\zeta)=\myt{f}_1(\zeta)$ for $\zeta\in D$.
An element $\myt{f}_N=(\zeta_N,\DD(\zeta_N,R_N))$ is called the {\it analytic continuation
of an element $\myt{f}_0$ along the path $\Gamma\subset\CC$} (in particular,
along the interval  $[0,6]$) if there exists a~{\it chain of elements}
$\myt{f}_j=(\zeta_j,\DD(\zeta_j,R_j))$, $j=1,\dots,N$, such that $\zeta_j\in\Gamma$ for all~$j$
and each $\myt{f}_{j+1}$ is a~direct analytic continuation of the element
$\myt{f}_j$, $j=0,\dots,N-1$. Clearly, in our examples
the Weierstrass extension from the point $x=0$ to the point $x=6$ along the positive real half-line
is possible because all four branch points of the function $\myt{f}$ lie away from the real line~$\RR$.

\subsection{}\label{s5s2}
We now proceed with numerically obtained  Examples~\ref{ex1}--\ref{ex3}.

\begin{figure}[!ht]
\centerline{
\includegraphics[width=15cm,height=15cm]{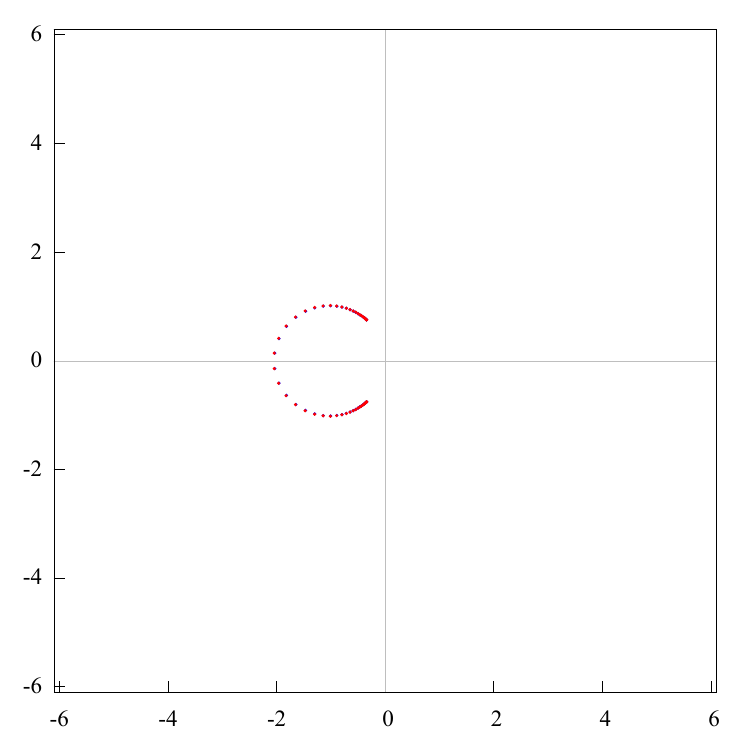}}
\vskip-6mm
\caption{Example~\ref{ex1}.
The zeros and poles (dark blue and red points) of the diagonal
PA $[50/50]_{\myt{f}_0}$, as constructed from the element
$\myt{f}_0\in\HH(0)$ from Example~\ref{ex1}. In this case, the cut
corresponding to the Stahl compact set lies entirely in the left half-plane. As a~result,
the value of the Weierstrass extension of the element $\myt{f}_0$ at the point $x=6$ can be
approximately evaluated using the PA $[50/50]_{\myt{f}_0}(6)$.
}
\label{fig1}
\end{figure}

\begin{example}\label{ex1}
In this example, the value of the parameter~$a$ in representation \eqref{62} is negative,
$a<0$. Correspondingly, both complex conjugate branch points
$a\pm ib$ lie in the left half-plane. Figure~\ref{fig1}
depicts the zeros and poles (dark blue and red points) of the PA $[50/50]_{\myt{f}_0}$, as constructed
from an element $\myt{f}_0\in\HH(0)$. In this case, the cut corresponding to the Stahl compact set lies
entirely in the left half-plane. Thus, the values of the Weierstrass extension
of the element $\myt{f}_0$ at the point $x=6$ can be approximately evaluated
using the PA $[50/50]_{\myt{f}_0}(6)$. Increasing the order of the PA  $[n/n]_{\myt{f}_0}$,
we will increase the accuracy of calculation of the value $\myt{f}_0(6)$ in accordance with~\eqref{8} and following well-known rules for
evaluating PA.
\end{example}

\begin{figure}[!ht]
\centerline{
\includegraphics[width=15cm,height=15cm]{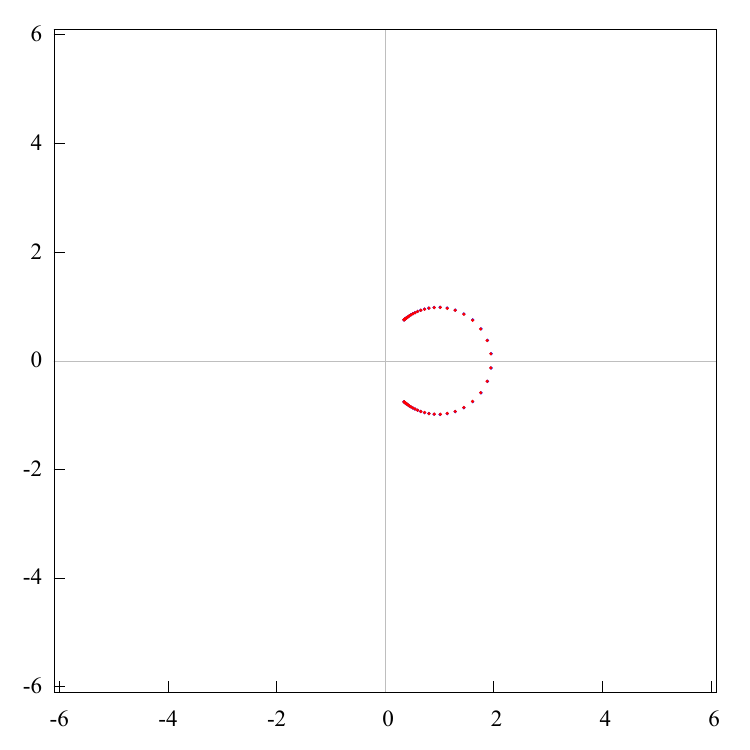}}
\vskip-6mm
\caption{Example~\ref{ex2}.
The zeros and poles (dark blue and red points) of the diagonal
PA $[50/50]_{\myt{f}_0}$, as constructed from the element
$\myt{f}_0\in\HH(0)$ from Example~\ref{ex2}. In this case, the cut corresponding
to the  Stahl compact set intersects the real line between the points $x=0$ and
$x=6$. Hence, the value
of the multivalued function~$\myt{f}$ at the point $x=6$, as
calculated from the PA $[n/n]_{\myt{f}_0}$ as $n\to\infty$, is different from the value
of the Weierstrass extension of the power series $\myt{f}_0$ from the point $x=0$ to the point $x=6$.
}
\label{fig2.1}
\end{figure}

\begin{figure}[!ht]
\centerline{
\includegraphics[width=15cm,height=15cm]{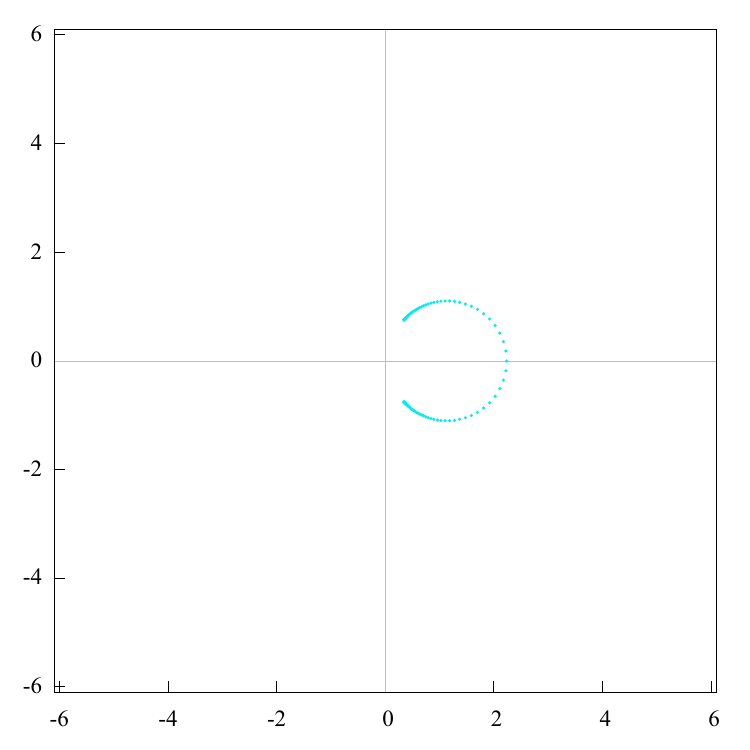}}
\vskip-6mm
\caption{Example~\ref{ex2}.
The zeros (light blue  points) of the type~II Hermite--Pad\'e polynomial
$Q_{2n}$ for $n=50$. By~\eqref{83}, the corresponding cut is the compact set
$E_2$ modeled by 100 zeros of this polynomial. This compact set $E_2$, as well as the
Stahl compact set, intersects the real line between the points $x=0$ and
$x=6$. Thus, from representation~\eqref{83} one cannot find the sought-for
value $\myt{f}_0(6)$ of the function.
But from the above it follows that from~\eqref{83} one can evaluate $\myt{f}_0(x^{(1)})$ at $x=6$
(we note once again that this value is distinct from the sought-for Weierstrass  value $\myt{f}_0(6)$).
}
\label{fig2.2}
\end{figure}

\begin{figure}[!ht]
\centerline{
\includegraphics[width=15cm,height=15cm]{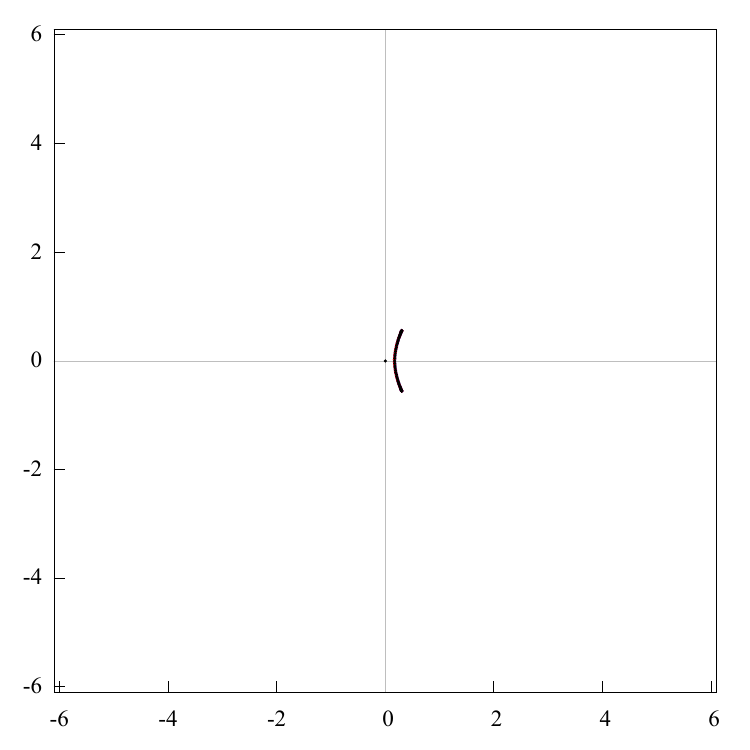}}
\vskip-6mm
\caption{Example~\ref{ex2}.
The zeros (red, dark blue and black points) of the type~I Hermite--Pad\'e polynomial $Q_{\nn,j}$ for $\nn=(50,50)$.
}
\label{fig2.3}
\end{figure}

\begin{figure}[!ht]
\centerline{
\includegraphics[width=15cm,height=15cm]{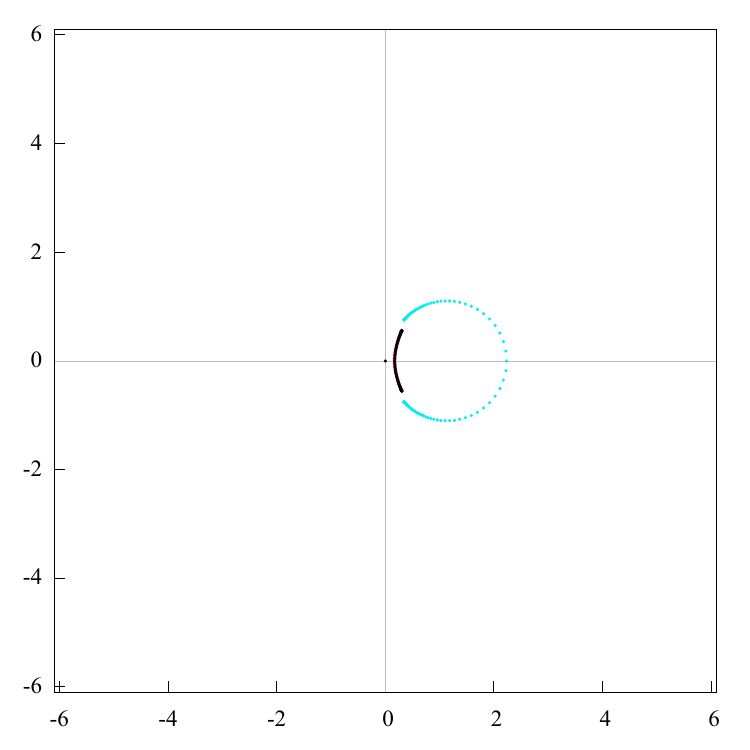}}
\vskip-6mm
\caption{Example~\ref{ex2}.
The figure combines Figs.~\ref{fig2.2} and~\ref{fig2.3}. It is seen that
even though the compact set $F_2$, which attracts the zeros of
the type~I Hermite--Pad\'e polynomials, intersects the real line between the points $x=0$ and $x=6$,
it still lies to the left of the compact set $E_2$.
It follows that the sought-for Weierstrass value $\myt{f}_0(6)$
agrees with the value $\myt{f}(x^{(2)})$ for $x=6$. Hence this value can be
evaluated by successive applications of relations~\eqref{83} and~\eqref{84}.
}
\label{fig2.4}
\end{figure}

\begin{example}\label{ex2}
In this example we choose a~positive $a$ in representation \eqref{62}. Thus, the pair of branch points $a\pm ib$
of the element $\myt{f}_0$ now lies in the right half-plane.
Evaluation of zeros and poles of the diagonal  PA $[50/50]_{\myt{f}_0}$ of the germ $\myt{f}_0$
shows (see Fig.~\ref{fig2.1}) that the cut corresponding
to the Stahl compact set intersects the real line between the points $x=0$ and
$x=6$. Therefore, the value of the
multivalued function~$f$ at the point $x=6$, as obtained from the PA $[n/n]_{\myt{f}_0}$ as $n\to\infty$,
does not agree with the value of the Weierstrass extension of the  power series $\myt{f}_0$ from the point $x=0$ to the point $x=6$.

Let us now employ the results of \S\,\ref{s4s3}.
We employ in our setting relations \eqref{83} and~\eqref{84} and
the above results on the limit distribution of the zeros of Hermite--Pad\'e polynomials.
They will be used to study the
topology of the Nuttall partition into sheets of that three-sheeted R.s.\ $\RS_3(w)$
which is associated to the germ  $\myt{f}_0$ of the function \eqref{62} (such a~surface is called a~Nuttall surface).
As was pointed out in~\S\,\ref{s4s3}, relations \eqref{83} and~\eqref{84} were obtained for a~germ given at
the point at infinity. For a~germ defined at the point $\zeta=0$, all these results
remain valid under the corresponding modification of the Nuttall partition into sheets
to the case under consideration of a~germ defined at the point $\zeta=0$.

Using this example, we shall show how the Hermite--Pad\'e polynomials, as given by \eqref{81} and~\eqref{82} and
relations \eqref{83} and~\eqref{84} can be used in the case when the diagonal PA are unfit for
numerical realization of the procedure for Weierstrass extension of a~given germ.

We cannot in fact study the structure of the boundaries
$F^{(1,2)}$ and $F^{(2,3)}$ between three sheets of the R.s.\ $\RS_3(w)$ directly with the help of Hermite--Pad\'e polynomials. However,
based on \eqref{83} and~\eqref{84} we can understand the structure of the projections
$E_2=\pi_3(F^{(1,2)})$ and $F_2=\pi_3(F^{(2,3)})$ of these boundaries to the Riemann
sphere, derive necessary information from this knowledge, and make appropriate conclusions on this basis.

Figure~\ref{fig2.2} shows the zeros (light blue points) of the type~II
Hermite--Pad\'e polynomial $Q_{2n}$ for $n=50$. By~\eqref{83}, the
corresponding cut is the compact set $E_2$ modeled by 100 zeros of this
polynomial. This compact set $E_2$, as well as the Stahl compact set,
intersects the real line between the points $x=0$ and $x=6$. This being
so, using \eqref{83} one cannot find the required value of the function
$\myt{f}_0(6)$. At first sight, the situation is similar to that
encountered in an attempt to employ the diagonal PA to find the value
$\myt{f}_0(6)$. However this is not so. From the above it follows that
with the help of \eqref{83} we can find the value $\myt{f}_0(x^{(1)})$
at $x=6$ (we note once again that this value does not agree with the
sought-for Weierstrass value $\myt{f}_0(6)$).

Let us now proceed to the type~I Hermite--Pad\'e polynomials $Q_{\nn,j}$ with $\nn=(50,50)$.
Their zeros  (red, dark blue, and black points) are depicted in Fig.~\ref{fig2.3}.
Figure~\ref{fig2.4} combines Figs.~\ref{fig2.2} and~\ref{fig2.3}. For the compact set $F_2$,
which attracts the zeros of type~I Hermite--Pad\'e polynomials,
Fig.~\ref{fig2.4} shows that even though
this set intersects the real line between the points $x=0$ and $x=6$,
but it still lies to the left of the compact set $E_2$. These compact sets $E_2$ and $F_2$, which are
projections of the boundaries between the sheets of the Nuttall surface $\RS_3(w)$, reflect the
following properties of the Nuttall partition of the  R.s.\ $\RS_3(w)$ into sheets from the
viewpoint of the Weierstrass extension of the original germ  $\myt{f}_0$.
When realizing the Weierstrass extension of the germ  $\myt{f}_0$ from the point $x=0$
to the point $x=6$ along the positive real half-line, we will cross at some point
$x=x_1\in(0,6)$ the compact set $E_2$. Thus, at this instant a~{\it transition from the
first sheet $\RS^{(1)}_3$ of the R.s.\ $\RS_3(w)$ to the second sheet $\RS^{(2)}_3$}
of this R.s.\ takes place. Further motion in accordance with the Weierstrass procedure to the point
$x=6$ proceeds already along the {\it second} sheet $\RS^{(2)}_3$ of the R.s.\ $\RS_3(w)$.
The compact set $F_2=\pi_3(F^{(2,3)})$ lies to the left of the compact set $E_2=\pi_3(F^{(1,2)})$.
Hence, when moving further to the right along the
second sheet $\RS^{(2)}_3$ to the point $x^{(2)}$, where $x=6$, we will always be only on the second sheet and
will not cross the boundary $F^{(2,3)}$ between the second and third sheets of the R.s.\ $\RS_3(w)$.
As a~result, the sought-for Weierstrass extension $\myt{f}_0(6)$
coincides with the value $\myt{f}(x^{(2)})$ at $x=6$. Consequently, this value can be
evaluated by a~repeated application of relations~\eqref{83} and~\eqref{84}.
\end{example}

\begin{example}\label{ex3}
Let us now proceed with the third example. In this example,
the parameters $a$ and~$b$ in representation~\eqref{62} are the same as
in Example~\ref{ex2}. However, the parameter $A_1$ is now altered. As in
Example~\ref{ex2}, the Stahl compact set for the germ $\myt{f}_0$ is the same
circular arc passing through the three points $\zeta=0,1/(a\pm ib)$.
As before, this arc intersects the positive real half-line  between the points $x=0$ and $x=6$.

Proceeding as in Example~\ref{ex2}, let us find for the germ $\myt{f}_0$
the zeros of the type~II Hermite--Pad\'e polynomial $Q_{2n}$ for $n=50$ and
the zeros of the type~I Hermite--Pad\'e polynomials $Q_{\nn,j}$, $j=0,1,2$, for $\nn=(50,50)$.
Using these data and employing relations \eqref{83} and~\eqref{84}, we analyze the
structure of the Nuttall partition into sheets associated with the germ
$\myt{f}_0$ of the three-sheeted R.s.\ $\RS_3(w)$. In accordance to~\eqref{83},
the zeros of the polynomial $Q_{2n}$ model the compact set $E_2=\pi_3(F^{(1,2)})$, while
the zeros of the polynomials
$Q_{\nn,j}$ model the compact set $F_2=\pi_3(F^{(2,3)})$. As in Example~\ref{ex2}, both these compact sets
intersect the positive real half-line
between the points $x=0$ and $x=6$. But now, in contrast to Example~\ref{ex2},
the compact set $F_2$ lies to the right of the compact set $E_2$ (see
Fig.~\ref{fig3.2}).
For the problem under consideration this means the following.
In the realization of the Weierstrass extension, when moving to the right of the point
$x=0$ along the real line to the point $x=6$, we meet the compact set $E_2$
at some\footnote{This
point $x_1$ is distinct from the point $x_1$ from Example~\ref{ex2}.}
$x=x_1$, and hence, the subsequent
realization of the Weierstrass extension will in fact proceed on the second sheet $\RS^{(2)}_3$
of the Nuttall surface $\RS_3(w)$.
As the motion continues further to the right from the point $x=x_1$ to the point $x=6$, we
intersect the compact set $F_2$ at some point $x=x_2$. This means that at this time we moved
from the second sheet of the R.s.\ $\RS_3(w)$ to its third sheet $\RS^{(3)}_3$.
Further realization of the Weierstrass procedure means that the actual motion now proceeds
already along the third sheet $\RS^{(3)}_3$
of the R.s.\ $\RS_3(w)$. However, as was already pointed out above,
the original germ
$\myt{f}_0$ cannot be continued on this third sheet as a~singlevalued analytic function,
because for such a~continuation appropriate cuts should be drawn on the third sheet.
Due to this obstacle, which is inherent in the realization of the Weierstrass procedure
by the compact set $F^{(2,3)}$ (the boundary between the second and third sheets of R.s.\ $\RS_3(w)$), we
cannot find the required value of $\myt{f}_0(6)$
using the Hermite--Pad\'e polynomials $Q_{2n}$ and $Q_{\nn,j}$ for $\nn=(n,n)$ in
the way we did it in Example~\ref{ex2},
because on the third sheet $\RS^{(3)}_3$ of the R.s.\ $\RS_3(w)$ relations \eqref{83}
and~\eqref{84} do not apply anymore.

Clearly, in this situation to realize the Weierstrass extension
one needs to augment the family of germs $[1,\myt{f}_0,\myt{f}_0^2]$ with the germ
$\myt{f}_0^3$,
increase the dimension of multiindices by~1, and proceed with constructions
\eqref{11}--\eqref{13},~\eqref{18} and~\eqref{24} corresponding to the
multiindices $\nn=(n,n,n)$, $\nn_1=(n,n-1,n-1)$, $\nn_2=(n,n,n-1)$ and
$\nn_3=(n,n,n)$. As a~result, on the four-sheeted R.s.\ $\RS_4(\myt{f})$ there appears
the Nuttall partition into sheets corresponding to the point $\zeta=0$ at which the original germ $\myt{f}_0$
is given.
The results of the corresponding numerical calculations with $n=50$ are shown in
Fig.~\ref{fig3.2}--\ref{fig3.5} (the zeros and poles of the diagonal PA
$[50/50]_{\myt{f}_0}$ are located quite similarly to Example~\ref{ex2},
and hence are not shown here).

Figure~\ref{fig3.3} depicts the zeros of the type~II Hermite--Pad\'e
polynomial $Q_{3n}$ for $n=50$ (yellow points), which model in
accordance with Theorem~\ref{th6} the  compact set
$E=\pi_4(\Gamma^{(1,2)})$, where $\Gamma^{(1,2)}$ is the boundary
between the first and second sheets of the R.s.\ $\RS_4(\myt{f})$. Besides,
Figure~\ref{fig3.3} shows (light blue points) the zeros of the
polynomial $S_{2n,1}$ for $n=50$, modeling the compact set
$F=\pi_4(\Gamma^{(2,3)})$, where $\Gamma^{(2,3)}$ is the boundary
between the second and third sheets of the R.s.\ $\RS_4(\myt{f})$. Figure
\ref{fig3.4} depicts  (black points) the zeros of the polynomial
$Q_{\nn,3}$, $\nn=(n,n,n)$, for $n=50$, which model the compact set
$E'=\pi_4(\Gamma^{(3,4)})$, where $\Gamma^{(3,4)}$ is the boundary
between the third and fourth sheets of the R.s.\ $\RS_4(\myt{f})$.
Figure~\ref{fig3.5} combines Figs.~\ref{fig3.3} and~\ref{fig3.4}. In
this case, the realization of the Weierstrass extension  from the point
$x=0$ to the point $x=6$ along the positive real half-line means in
fact the following. Moving to the right from the point $x=0$ to the
point $x=6$ we meet the compact set~$E$ at some point\footnote{This
point $x_1$ and the point $x_2$ are in general different from those
points $x_1$ and $x_2$ used above in the present paper and in
Example~\ref{ex2}.} $x=x_1$. This means that when realizing the
Weierstrass extension our further motion proceeds now along the second
sheet $\RS^{(2)}_4$ of the R.s.\ $\RS_4(\myt{f})$. This continues until at
some point $x=x_2$ we intersect the compact set~$F$, which is the
projection of the boundary $\Gamma^{(2,3)}$ between the second and
third sheets of the R.s.\ $\RS_4(\myt{f})$. From this point on, our motion to
the right takes place already along the third sheet $\RS^{(3)}_4$ of
the R.s.\ $\RS_4(\myt{f})$. The fact that the compact set~$E'$, which is the
projection of the boundary between the third and fourth sheets of the
R.s.\ $\RS_4(\myt{f})$, lies to the left of the  compact set~$F$ (see
Fig.~\ref{fig3.5}, which combines Figs.~\ref{fig3.3} and~\ref{fig3.4}),
means for us that when moving to the right along the third sheet
$\RS^{(3)}_4$ of the R.s.\ $\RS_4(\myt{f})$ we shall not meet the boundary
between the third and fourth $\RS^{(4)}_4$ sheets of the R.s.\
$\RS_4(\myt{f})$. So, during the remaining time (until the point $x=6$ is
reached) we shall in fact move along the third sheet. Note that the
compact set $E'$, as well as the compact sets $E$ and~$F$, intersects
the real line between the points $x=0$ and $x=6$ and lies to the left
of the compact set~$E$. This, however, has no effect on the
implementation of the Weierstrass procedure, because the compact set
$E'$ is the projection of the boundary $\Gamma^{(3,4)}$ between the
third and fourth sheets, which ``is not seen'' on the first and second
sheets of the R.s.\ $\RS_4(\myt{f})$.

Thus, according to what has been said, the required Weierstrass value
$\myt{f}_0(6)$ coincides with the value $\myt{f}(x^{(3)})$ of the function at $x=6$.
According to the above Theorems~\ref{th5}--\ref{th7}, the required value of
$\myt{f}_0(6)$ is evaluated by successive application of
relations\footnote{Recall that all Hermite--Pad\'e polynomials involved in these relations
should be calculated from the germ $\myt{f}_0$ of the function
$\myt{f}$ of the form \eqref{44} at the point $\zeta=0$.} \eqref{56},~\eqref{48.2} and
~\eqref{47} (see also~\eqref{91}--\eqref{93}).
\end{example}

\begin{remark}\label{rem8}
To conclude we note that the new method proposed here of
efficient continuation of a~power series beyond its disc of convergence
guarantees the possibility, in the class of
multivalued analytic functions with finite number of branch points,
to realize the Weierstrass procedure arbitrarily far\footnote{Of course,
if one is concerned with a~specific continuation, for example,
along the positive real half-line, then such a~continuation is possible only up to the nearest
branch point.} beyond the disc of convergence. To this end one needs to
appropriately increase the dimensions of the employed multiindices.

By now there are practically no theoretical results justifying the use of the
method of analytic continuation proposed here.
We plan to offer the proofs of the above Theorems~\ref{th1}--\ref{th3} and
Theorems~\ref{th5}--\ref{th7} in the separate papers~\cite{IkKoSu18} and \cite{LoKoSu18}, respectively.
We propose to carry out numerical analysis of applied problems on the
basis of the method proposed in the present paper.
\end{remark}

It is worth pointing out that the new approach to the solution of the problem of
efficient analytic continuation of power series was inspired not only by theoretical results and
conjectures mentioned above, but also by the numerical experiments, whose results are given in the papers
\cite{IkKoSu15},~\cite{IkKoSu15b} and~\cite{IkKoSu16}.



\clearpage
\newpage

\clearpage
\newpage

\clearpage
\newpage

\clearpage
\newpage

\clearpage
\newpage

\clearpage
\newpage

\clearpage
\newpage

\clearpage
\newpage

\clearpage
\newpage
\begin{figure}[!ht]
\centerline{
\includegraphics[width=15cm,height=15cm]{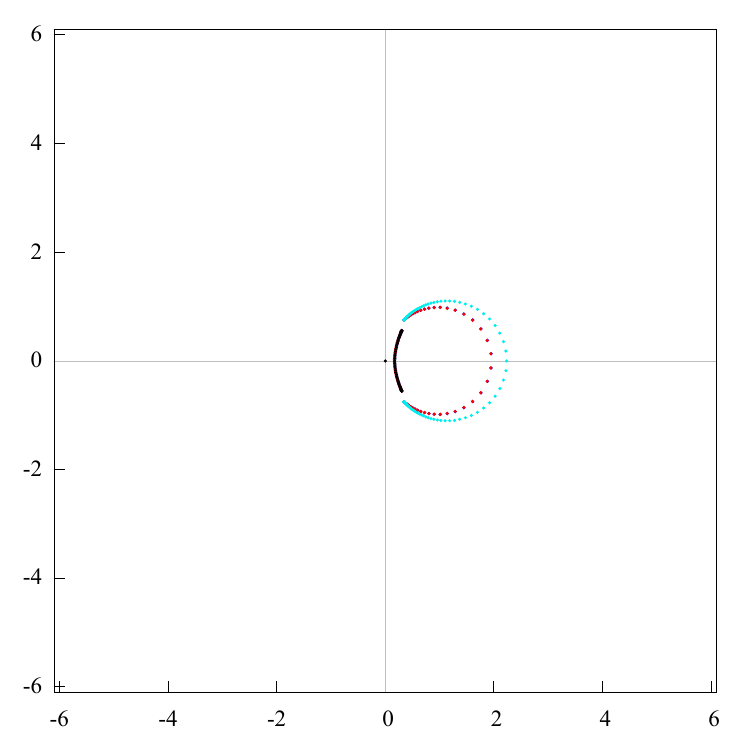}}
\vskip-6mm
\caption{Example~\ref{ex2}.
The zeros and poles (dark blue and red points) of the
PA $[50/50]_{\myt{f}_0}$, as constructed from the element
$\myt{f}_0\in\HH(0)$ from Example~\ref{ex2}. The figure also shows the zeros
(light blue  points) of the type~II Hermite--Pad\'e polynomial
$Q_{2n}$ for $n=50$ and the zeros (red, dark blue and black points) of the type~I Hermite--Pad\'e polynomials
$Q_{\nn,j}$ for $\nn=(50,50)$.
It is seen that all three compact sets $S$, $E_2$, and $F_2$ are distinct from each other.
}
\label{fig2.5}
\end{figure}


\clearpage
\newpage
\begin{figure}[!ht]
\centerline{
\includegraphics[width=15cm,height=15cm]{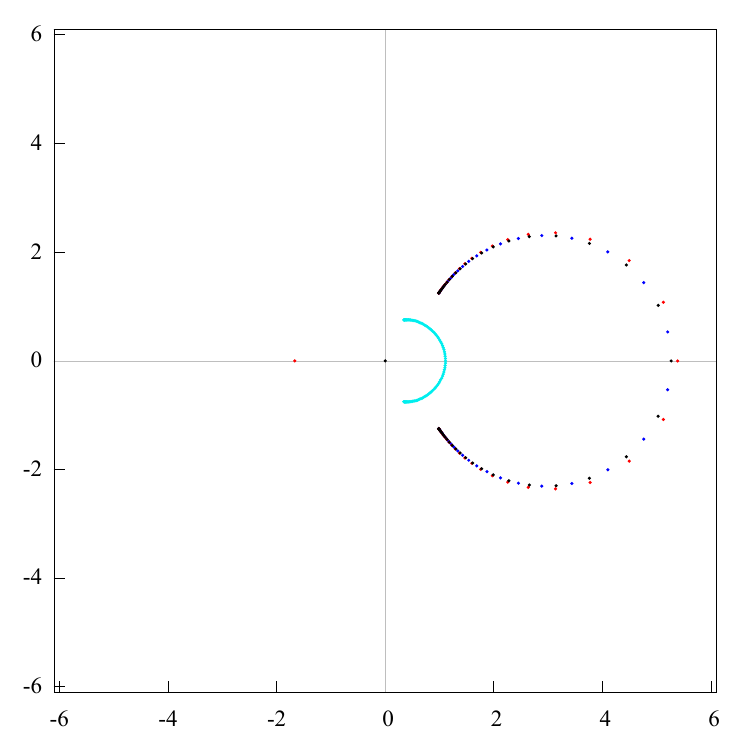}}
\vskip-6mm
\caption{Example~\ref{ex3}.
The zeros of the type~II Hermite--Pad\'e polynomial $Q_{2n}$ (light blue  points)
for $n=50$ and the zeros (red, dark blue and black points)  of the type~I Hermite--Pad\'e polynomials $Q_{\nn,j}$, $j=0,1,2$,
for $\nn=(50,50)$.
From these data and with the help of relations \eqref{83} and~\eqref{84}
the structure of the Nuttall partition into sheets of
the three-sheeted R.s.\ $\RS_3(w)$ with the germ $\myt{f}_0$ is analyzed.
In accordance with \eqref{83}, the zeros of the polynomial $Q_{2n}$ model the
compact set $E_2=\pi_3(F^{(1,2)})$ and the zeros of the polynomials
$Q_{\nn,j}$ model the compact set $F_2=\pi_3(F^{(2,3)})$. Both these compact sets
intersect the real line, but, in contrast to Example~\ref{ex2}
(Fig.~\ref{fig2.4}), the compact set $F_2$ lies to the right of the compact  set~$E_2$.
}
\label{fig3.2}
\end{figure}

\clearpage
\newpage
\begin{figure}[!ht]
\centerline{
\includegraphics[width=15cm,height=15cm]{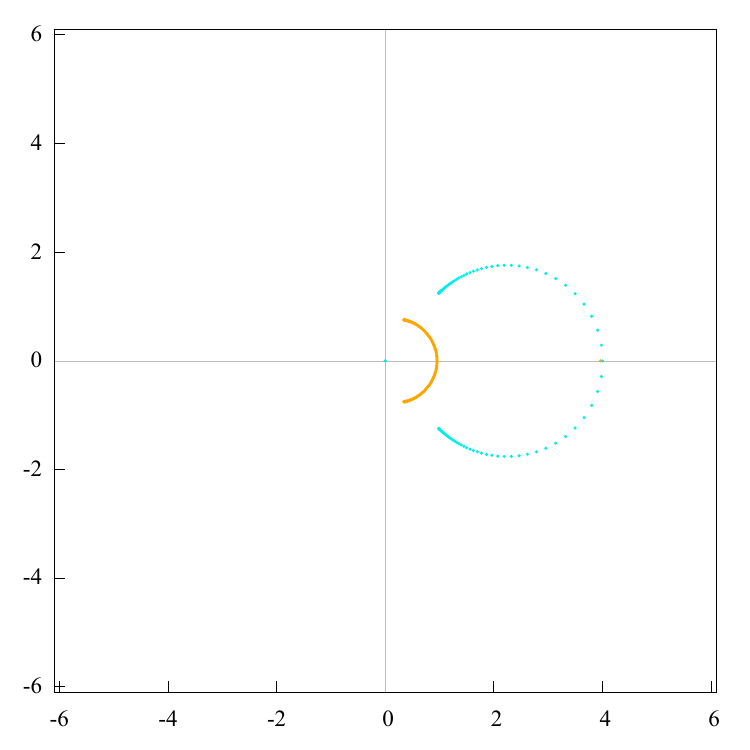}}
\vskip-6mm
\caption{Example~\ref{ex3}.
The zeros (yellow points) of the type~II Hermite--Pad\'e polynomial
$Q_{3n}$ for $n=50$, which model in accordance with
Theorem~\ref{th6} the compact set $E=\pi_4(\Gamma^{(1,2)})$.
Moreover, the figure depicts the zeros  (light blue  points) of the polynomials  $S_{2n,1}$ for $n=50$, which model the compact set
$F=\pi_4(\Gamma^{(2,3)})$.
}
\label{fig3.3}
\end{figure}

\clearpage
\newpage
\begin{figure}[!ht]
\centerline{
\includegraphics[width=15cm,height=15cm]{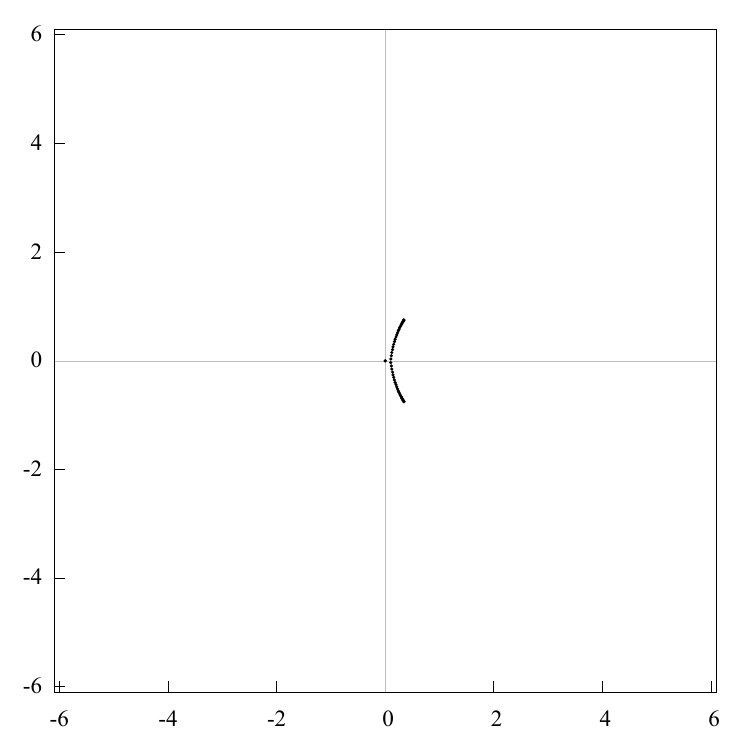}}
\vskip-6mm
\caption{Example~\ref{ex3}.
The zeros (black points) of the polynomial $Q_{\nn,2}$, $\nn=(n,n,n)$, for
$n=50$, which model the compact set
$E'=\pi_4(\Gamma^{(3,4)})$.
}
\label{fig3.4}
\end{figure}

\clearpage
\newpage
\begin{figure}[!ht]
\centerline{
\includegraphics[width=15cm,height=15cm]{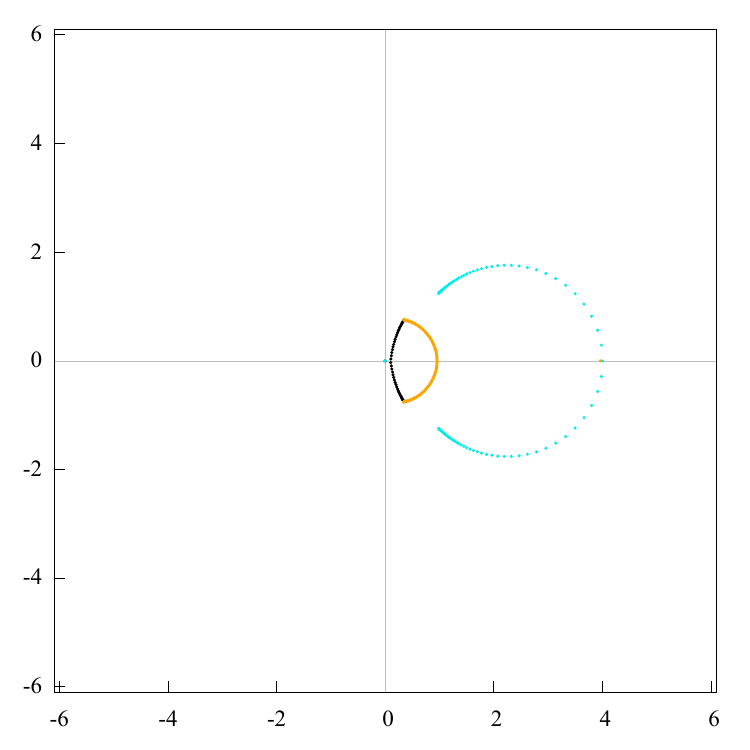}}
\vskip-6mm
\caption{Example~\ref{ex3}.
The figure combines Figs.~\ref{fig3.3} and~\ref{fig3.4}.
It shows that all three compact sets $E$, $F$ and $E'$ intersect the real line between the points $x=0$
and $x=6$. But the compact set $E'$ lies to the left of the compact sets $E$ and $F$.
}
\label{fig3.5}
\end{figure}

\clearpage
\newpage
\begin{figure}[!ht]
\centerline{
\includegraphics[width=15cm,height=15cm]{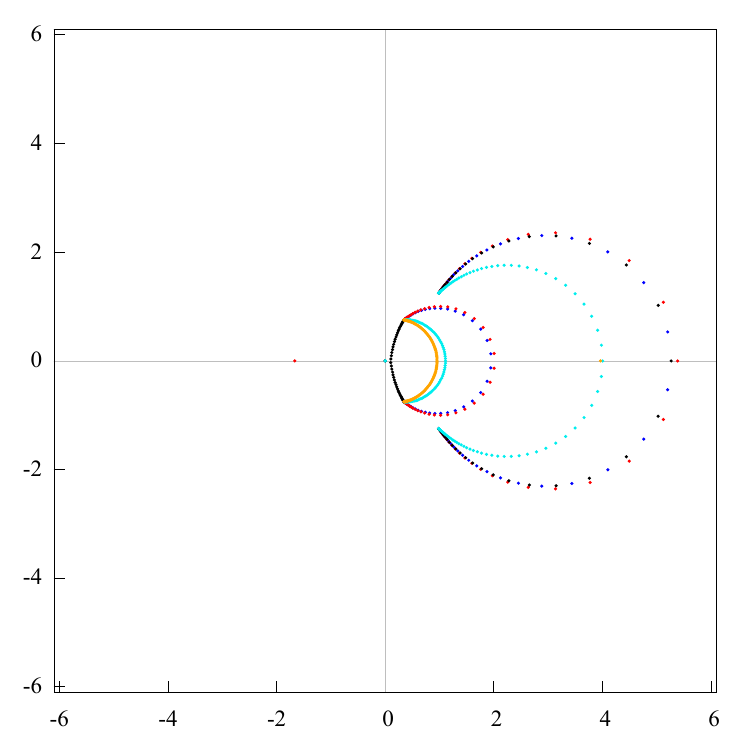}}
\vskip-6mm
\caption{Example~\ref{ex3}.
The zeros and poles (dark blue and red points)
PA $[50/50]_{\myt{f}_0}$, as constructed from the element
$\myt{f}_0\in\HH(0)$ from Example~\ref{ex3}. Moreover, the figure depicts
the zeros (light blue  points) of the type~II Hermite--Pad\'e polynomial
$Q_{2n}$ for $n=50$ and the zeros (red, dark blue and black points) of the
type~I Hermite--Pad\'e polynomials $Q_{\nn,j}$ for $\nn=(50,50)$ (Fig.~\ref{fig3.2}).
Besides, it shows the zeros (yellow points) of the type~II Hermite--Pad\'e polynomial
$Q_{3n}$, the zeros (yellow points) of the polynomial $S_{2n,1}$, and the zeros (black points)
of the type~I Hermite--Pad\'e polynomial $Q_{\nn,2}$ for $\nn=(n,n,n)$, where $n=50$.
It is clear that all six compact sets $S$, $E_2,F_2$, $E,F$ and $E'$ are distinct from each other.
}
\label{fig3.6}
\end{figure}


\end{document}